\RequirePackage{fix-cm}
\documentclass[smallextended]{svjour3}
\usepackage{amssymb}
\usepackage{amsfonts}
\usepackage{amsmath}
\usepackage{theorem}
\usepackage{pgfplots}
\usepackage{hyperref} 
\hypersetup{plainpages=false, colorlinks, linkcolor=black, citecolor=black, urlcolor=blue,
pdftitle={NFFT error},
pdfauthor={Daniel Potts, Manfred Tasche},
pdfstartview={FitBH}}

\theoremstyle{plain}
\newtheorem{algorithm}[theorem]{Algorithm}

\newcommand{\bend}{\hspace*{0ex} \hfill \hbox{\vrule height
    1.5ex\vbox{\hrule width 1.4ex \vskip 1.4ex\hrule  width 1.4ex}\vrule
    height 1.5ex}}
\newcommand{\qedsymbol}{\rule{1.5ex}{1.5ex}}

\newenvironment{Lemma}{\goodbreak\begin{lemma}\rmfamily}{\end{lemma}}
\newenvironment{Theorem}{\goodbreak\begin{theorem}}{\end{theorem}}
\newenvironment{Remark}{\goodbreak\begin{remark}\rmfamily}{\bend\end{remark}}

\newenvironment{Example}{\goodbreak\begin{example}\rmfamily}{\bend\end{example}}

\numberwithin{equation}{section}
\numberwithin{table}{section}
\numberwithin{figure}{section}

\renewcommand{\mathbf}[1]{\ensuremath{\boldsymbol{#1}}}

\journalname{Preprint TU Chemnitz}
\begin{document}

\title{Uniform error estimates for nonequispaced fast Fourier transforms}

\author{Daniel Potts \and Manfred Tasche}

\authorrunning{D.~Potts, M.~Tasche} %

\institute{D. Potts \at
	Chemnitz University of Technology, Faculty of Mathematics, 09107 Chemnitz, Germany \\
	\email{potts@mathematik.tu-chemnitz.de}
	\and
	M. Tasche \at
	University of Rostock, Institute of Mathematics, D--18051 Rostock, Germany\\
	\email{manfred.tasche@uni-rostock.de}
}

\date{Received: date / Accepted: date}

\maketitle

\begin{abstract}
In this paper, we study the error behavior of the nonequispaced fast Fourier transform (NFFT). This approximate algorithm is mainly based on the convenient choice of a compactly supported window function.
So far, various window functions have been used and new window functions have recently been proposed. We present novel error estimates for NFFT with compactly supported, continuous window functions and derive rules for convenient choice from the parameters involved in NFFT. The error constant of a window function depends mainly on the oversampling factor and the truncation parameter.

\medskip

\emph{Key words}: nonequispaced fast Fourier transform, NFFT, error estimate, Wiener algebra, oversampling factor, truncation parameter,
compactly supported window function, Paley--Wiener theorem.
\smallskip

AMS \emph{Subject Classifications}:
65T50,  
94A12,  
42A10.   
\end{abstract}


\section{Introduction}
Since the restriction to equispaced data is an essential drawback in several applications of discrete Fourier transform, one has developed fast algorithms for nonequispaced data, the so-called \emph{nonequispaced fast Fourier transform} (NFFT), see \cite{DuRo93,Bey95,St97,postta01,KeKuPo09,nfft3} and \cite[Chapter 7]{PlPoStTa18}.

In this paper we investigate error estimates for the NFFT, where we restrict ourselves
to an approximation by translates of a previously selected, continuous window function with compact support. Other approaches based on the fast multipole method and on the low-rank approximation were presented in \cite{duro95,RuTo18}.

After the seminal paper \cite{St97}, the similarities of the window-based algorithms for NFFT became clear. In the following, we give an overview of the window-based NFFT used so far. In this construction of NFFT, a window function is applied together with its Fourier transform.  This connection is very important in order to deduce error estimates for the NFFT. This made it possible to determine convenient parameters of the involved window function.  To develop an NFFT, the necessary Fourier coefficients of the periodized window function can simply be calculated by a convenient quadrature rule. A challenge are error estimations in order to determine the parameters of the
window function involved. By $C(\mathbb T)$ we denote the Banach space of all 1-periodic, continuous functions, where $\mathbb T = {\mathbb R}/{\mathbb Z}$ is the torus.

{The considered window functions depend on some parameters. Assume that $N \in 2 \mathbb N$ is the order of the given 1-periodic trigonometric polynomial which values will be computed by NFFT. Let $\sigma > 1$ be an \emph{oversampling factor} such that $N_1:= \sigma N \in 2 \,\mathbb N$. For fixed \emph{truncation parameter} $m \in {\mathbb N}\setminus \{1\}$
with $2m \ll N_1$,
we denote by $\Phi_{m,N_1}$ the set of all window functions $\varphi:\, \mathbb R \to [0,\,1]$ with the following properties:
\begin{quote}
$\bullet$ Each window function $\varphi$ is even, has a  compact support $[-m/N_1,\,m/N_1]$, and is continuous on $\mathbb R$.\\
$\bullet$ Each restricted window function $\varphi|_{ [0,\, m/N_1]}$ is decreasing with $\varphi(0) = 1$ and $\varphi(m/N_1) = 0$.\\
$\bullet$ For each window function $\varphi$, the Fourier transform
$$
{\hat \varphi}(v) := \int_{\mathbb R} \varphi(t)\, {\mathrm e}^{-2\pi {\mathrm i}\,v t}\, {\mathrm d}t = 2\, \int_0^{m/N_1}\varphi(t)\, \cos(2\pi\,v t)\, {\mathrm d}t
$$
is positive for all $v\in [-N/2,\,N/2]$.
\end{quote}
Note that for fixed $N_1$ the truncation parameter $m$ determines the size of the support of $\varphi \in \Phi_{m,N_1}$. If a window function $\varphi \in \Phi_{m,N_1}$ has the form
$$
\varphi(x) = \frac{1}{\varphi_1\big(\beta\,\varphi_2(0)\big)}\,\varphi_1\big(\beta\,\varphi_2(x)\big)\,, \quad x \in \mathbb R\,,
$$
with $\beta > 0$ and convenient functions $\varphi_1$, $\varphi_2$, then $\beta$ is a so-called \emph{shape parameter} of $\varphi$. Examples of window functions of the set $\Phi_{m,N_1}$ are the B-spline window function \eqref{eq:B-spline}, the modified B-spline window function \eqref{eq:modB-spline}, the algebraic window function \eqref{eq:algwindow}, the Bessel window function (5.14), the $\sinh$-type window function \eqref{eq:sinhwindow}, the related $\sinh$-type functions, see Subsection \ref{Sec:window}, the modified $\cosh$-type window function \eqref{eq:mcoshwindow} and the related $\cosh$-type window functions, see Subsection \ref{Sec:modwindow}. Note that the
Kaiser--Bessel window function (see \cite[p. 393]{PlPoStTa18}) and the Gaussian window function (see \cite[p. 390]{PlPoStTa18}) are not contained in $\Phi_{m,N_1}$, since these window functions
are supported on whole $\mathbb R$.
\medskip

The aim of this paper is a systematic approach to uniform error estimates for NFFT, where a compactly supported, continuous window function $\varphi \in \Phi_{m,N_1}$ is used. We introduce the $C(\mathbb T)$-error constant
$$
e_{\sigma}(\varphi) = \sup_{N \in 2 \mathbb N}\Big(\max_{n\in I_N}\big\| \sum_{r \in {\mathbb Z}\setminus \{0\}} \frac{{\hat \varphi}(n + rN_1)}{{\hat\varphi}(n)}\,{\mathrm e}^{2 \pi {\mathrm i}\,r N_1\, \cdot}\big\|_{C(\mathbb T}\Big)\,,
$$
where $I_N$ denotes the index set $\{-N/2,\,1-N/2,\,\ldots,\,N/2-1\}$.
As shown in Lemma \ref{lemma:convenient2}, the uniform error
of the NFFT with nonequispaced spatial data and equispaced frequencies can be estimated by $e_{\sigma}(\varphi)$.  Analogously in Lemma \ref{lemma:convenientNFFT}, the error of the NFFT
with nonequispaced frequencies and equispaced spatial data is estimated by $e_{\sigma}(\varphi)$ too. Therefore in the following, we study mainly the behavior of the $C(\mathbb T)$-error constant $e_{\sigma}(\varphi)$.
Our main result is Theorem \ref{Lemma:esigmaNvarphibounded}, where we describe a general concept for the construction of a convenient upper bound  for the $C(\mathbb T)$-error constant $e_{\sigma}(\varphi)$ with a window function $\varphi \in \Phi_{m,N_1}$.
Applying Theorem \ref{Lemma:esigmaNvarphibounded}, we obtain upper bounds for $e_{\sigma}(\varphi)$ with special window function $\varphi \in \Phi_{m,N_1}$.
We show that the $C(\mathbb T)$-error constant $e_{\sigma}(\varphi)$ of a window function $\varphi \in \Phi_{m,N_1}$ depends mainly on the oversampling factor $\sigma > 1$ and the truncation parameter $m \in \mathbb N \setminus \{1\}$. Since we are interested in NFFT with relatively low computational cost, the
oversampling factor $\sigma \in \big[\frac{5}{4},\,2\big]$ and the truncation parameter $m \in \{2,\,3,\, \ldots,6\}$ are restricted. These parameters $\sigma$ and $m$ determine the shape parameter $\beta$ of the window function. For the Bessel window function \eqref{eq:Besselwindow}, the $\sinh$-type window function \eqref{eq:sinhwindow}, and the modified $\cosh$-type window function \eqref{eq:mcoshwindow}, a good choice is the shape parameter $\beta = 2\pi m\, \big(1 - \frac{1}{2 \sigma}\big)$.}
\smallskip

In connection with NFFT, B-\emph{spline window functions} were first investigated in \cite{Bey95}.
In the important application of particle simulation (see \cite{DeHo98a}), the B-spline window function was also used. Later it became clear that these methods can be interpreted as a special case of the fast summation method, see \cite{post02,nestlerdiss} and the references therein. Based on this unified approach, one can use all the other window functions for this application too.
The convenient choice of the shape parameter is of special importance, as shown in \cite{Ne14} for the root mean square error of the NFFT.
\smallskip

In this paper we suggest four new continuous, compactly supported window functions for the NFFT, namely the algebraic, Bessel, $\sinh$-type, and modified $\cosh$-type window function. The \emph{algebraic window function} is very much related to the B-spline window function, but much simpler to compute. We show that the \emph{Bessel window function} (5.14), $\sinh$-\emph{type window function} (5.21), and \emph{modified} $cosh$-\emph{type window function} \eqref{eq:mcoshwindow} are very convenient for NFFT, since they possess very small $C(\mathbb T)$-error constants with exponential decay with respect to $m$.
{  It is difficult to design a window function $\varphi\in \Phi_{m,N_1}$ with minimal $C(\mathbb T)$-error constant.
We prove that the best error behavior has the modified $\cosh$-type window function \eqref{eq:mcoshwindow} with the shape parameter $\beta = 2 \pi m \big(1 - \frac{1}{2\sigma}\big)$, $\sigma \ge \frac{5}{4}$.}
Further we compare several window functions with respect to the corresponding $C(\mathbb T)$-error constants for the NFFT.  Based on the error estimates of {  the} $\sinh$-type window function, we are able to extend the error estimates, see \cite{PoTa20}, to window functions where an analytical expression of its Fourier transform is unknown, see \cite{FINUFFT,BaMaKl18,Ba20}.
\smallskip

{  We prefer the use of compactly supported, continuous window functions $\varphi \in \Phi_{m,N_1}$ by following reasons:
\begin{quote}
$\bullet$ As explained in Remark 1, the NFFT with a window function $\varphi \in \Phi_{m,N_1}$ is simpler than the NFFT with a window function supported on whole $\mathbb R$.\\
$\bullet$ The window functions $\varphi \in \Phi_{m,N_1}$ presented in Section 5 (and their Fourier transforms) have simple explicit forms and they are convenient as window functions for NFFT.\\
$\bullet$ Few window functions $\varphi \in \Phi_{m,N_1}$ (such as Bessel, $\sinh$-type, and modified $\cosh$-type window function) possess low $C(\mathbb T)$-error constants with exponential decay with respect to $m$.
The best error behavior has the modified $\cosh$-type window function.
\end{quote}
}
\medskip

The outline of the paper is as follows. In Section \ref{Sec:Window} we introduce the basic definitions and develop the error estimates for an NFFT with a general compactly supported, continuous window function.
Important tools for the estimation of the Fourier transforms of window functions are developed in Section \ref{Sec:FT}. In Section \ref{Sec:Paley-Wiener} we present a modified Paley--Wiener Theorem which characterizes the behavior of Fourier transforms {  of compactly supported functions} lying in a special Sobolev space.
 The main results of this paper are contained in Section \ref{Sec:SpWindow}. Using the uniform norm, we present explicit error estimates for
the (modified) B-spline, algebraic, Bessel, $\sinh$-type, and modified $\cosh$-type window functions. Further we show numerical tests so that the $C(\mathbb T)$-error constants of the different window functions can be easily compared.

\section{Convenient window functions for NFFT}\label{Sec:Window}

Let $\tilde{\varphi}:\, \mathbb T \to [0,\,1]$ be the 1-\emph{periodization} of $\varphi$, i.e.,
\begin{equation} \label{eq:tildevarphi}
\tilde{\varphi}(x) := \sum_{k\in \mathbb Z} \varphi(x+k)\,, \quad x \in \mathbb T\,.
\end{equation}
Note that for each $x\in \mathbb R$ the series \eqref{eq:tildevarphi} has at most one nonzero term.
Then the Fourier coefficients of $\tilde{\varphi}$ read as follows
$$
c_k(\tilde \varphi) := \int_0^1 \tilde{\varphi}(t) \,{\mathrm e}^{-2\pi{\mathrm i}\,k t}\,{\mathrm d}t =  {\hat \varphi}(k)\,, \quad k \in \mathbb Z\,.
$$
By the properties of the window function $\varphi \in \Phi_{m,N_1}$, the 1-periodic function $\tilde{\varphi}$ is continuous on $\mathbb T$ and of bounded variation over $\big[- \frac{1}{2},\, \frac{1}{2}\big]$. Then from
the Convergence Theorem of Dirichlet--Jordan (see \cite[Vol.~1, pp.~57--58]{Zy}), it follows that $\tilde{\varphi}$ possesses the uniformly convergent Fourier expansion
\begin{equation}
\label{eq:Fourierseriestildevarphi}
\tilde{\varphi}(x) = \sum_{k\in \mathbb Z} c_k(\tilde \varphi)\,{\mathrm e}^{2\pi{\mathrm i}\,k x}\,, \quad x \in \mathbb R\,.
\end{equation}

\begin{Lemma}\label{lemma2.1}
For given $\varphi \in \Phi_{m,N_1}$, the series
\begin{equation*}
\sum_{r \in {\mathbb Z}} c_{n+r N_1}(\tilde \varphi)\,{\mathrm e}^{2\pi {\mathrm i}\,(n + r N_1)\,x}
\end{equation*}
is convergent for each $x\in \mathbb R$ and has the sum
$$
\frac{1}{N_1}\, \sum_{\ell =0}^{N_1 -1} {\mathrm e}^{-2\pi {\mathrm i}\,n \ell/N_1}\, {\tilde \varphi}
(x + \frac{\ell}{N_1})
$$
which coincides with the rectangular rule of the integral
\begin{equation*}
c_n({\tilde \varphi}(x+ \cdot))=\int_0^1 {\tilde \varphi}(s + x)\, {\mathrm e}^{-2\pi {\mathrm i}\,ns}\,{\mathrm d}s = c_n(\tilde \varphi)\,{\mathrm e}^{2\pi {\mathrm i}\,nx}\,.
\end{equation*}
\end{Lemma}

\emph{Proof.} From \eqref{eq:Fourierseriestildevarphi} it follows that for all $n \in \mathbb Z$ and $x \in \mathbb R$ it holds
$$
{\mathrm e}^{-2\pi {\mathrm i}\,nx} \,{\tilde \varphi}(x) = \sum_{k\in \mathbb Z} c_{k+n}(\tilde \varphi)\,{\mathrm e}^{2\pi {\mathrm i}\,kx}\,.
$$
Replacing $x$ by $x + \frac{\ell}{N_1}$ with $\ell = 0,\ldots, N_1 -1$, we obtain
$$
{\mathrm e}^{-2\pi {\mathrm i}\,n\,(x + \ell/N_1)} \,{\tilde \varphi}(x + \frac{\ell}{N_1}) = \sum_{k\in \mathbb Z} c_{k+n}(\tilde \varphi)\,{\mathrm e}^{2\pi {\mathrm i}\,kx}\,{\mathrm e}^{2\pi {\mathrm i}\,k\ell/N_1}\,.
$$
Summing the above formulas for $\ell = 0,\ldots,N_1 - 1$ and using the known relation
$$
\sum_{\ell =0}^{N_1- 1} {\mathrm e}^{2\pi {\mathrm i}\,k\ell/N_1} = \left\{ \begin{array}{ll} N_1 & \quad k \equiv 0 \, {\mathrm{mod}}\,N_1\,,\\
0 & \quad k \not\equiv 0 \, {\mathrm{mod}}\,N_1\,,
\end{array} \right.
$$
we conclude that
$$
\sum_{\ell =0}^{N_1 - 1}{\mathrm e}^{-2\pi {\mathrm i}\,n\,(x + \ell/N_1)} \,{\tilde \varphi}(x + \frac{\ell}{N_1}) = N_1\,\sum_{r\in \mathbb Z} c_{n + r N_1}(\tilde \varphi)\,{\mathrm e}^{2\pi {\mathrm i}\,r N_1 x}\,.
$$
This completes the proof. \qedsymbol

\subsection{NFFT with noneqispaced spatial data and equispaced frequencies}

Let $M\in 2 \mathbb N$ and $I_M := \{-M/2, 1- M/2,\ldots ,M/2-1\}$ be given.
The NFFT \emph{with nonequispaced spatial data and equispaced frequencies} is an approximate, fast algorithm which computes approximately the values $f(x_j)$, $j \in I_M$,  of a 1-periodic trigonometric polynomial
\begin{equation}
\label{eq:f(x)}
f(x) := \sum_{k\in I_N} c_k(f)\, {\mathrm e}^{2\pi {\mathrm i}\,k x}
\end{equation}
at $M$ nonequispaced nodes $x_j \in [-\frac{1}{2},\,\frac{1}{2})$, $j \in I_M$. Using a window function $\varphi \in \Phi_{m,N_1}$, the trigonometric polynomial $f$
is approximated by the 1-periodic function
\begin{equation}
\label{eq:s(x)}
s(x) := \sum_{\ell \in I_{N_1}} g_{\ell}\, {\tilde \varphi}( x - \frac{\ell}{N_1})
\end{equation}
with conveniently chosen coefficients $g_{\ell}\in {\mathbb C}$. The computation of the values $s(x_j)$, $j\in I_M$, which approximate $f(x_j)$ is very easy. Since $\varphi$ is compactly supported and
thus $\tilde \varphi$ is well-localized, each value $s(x_j)$, $j\in I_M$, is equal to a sum of few nonzero terms.

The coefficients $g_{\ell}$ can be determined by discrete Fourier transform  (DFT) as follows.
The 1-periodic function $s$ possesses the Fourier expansion
$$
s(x) = \sum_{k\in \mathbb Z} c_k(s)\,{\mathrm e}^{2\pi {\mathrm i}\,kx}
$$
with the Fourier coefficients
$$
c_k(s) = \int_0^1 s(t)\, {\mathrm e}^{-2\pi {\mathrm i}\,k t}\, {\mathrm d}t = {\hat g}_k\,c_k({\tilde \varphi})\,, \quad k \in \mathbb Z\,,
$$
where
$$
{\hat g}_k := \sum_{\ell \in I_{N_1}} g_{\ell}\, {\mathrm e}^{- 2\pi {\mathrm i}\, k \ell/N_1}\,.
$$
In other words, the vector $({\hat g}_k )_{k \in I_{N_1}}$ is the DFT of length $N_1$ of coefficient vector $(g_{\ell})_{\ell \in I_{N_1}}$ such that
$$
{\hat g}_k  = {\hat g}_{k + N_1}\,, \quad k\in \mathbb Z\,.
$$
In order to approximate $f$ by $s$, we set
$$
{\hat g}_k = \left\{ \begin{array}{ll} \frac{c_k(f)}{c_k(\tilde \varphi)} & \quad k \in I_N\,,\\ [1ex]
0 & \quad k \in I_{N_1} \setminus I_N\,.
\end{array} \right.
$$
{  Note that the values ${\hat g}_k$, $k \in I_N$, can be used in an efficient way. Even if \eqref{eq:f(x)} is only known at finitely many equispaced points of $[0,\,1]$, the Fourier coefficients $c_k(f)$, $k \in I_N$,
can be approximately determined by a fast Fourier transform (FFT). For many window functions $\varphi \in \Phi_{m,N_1}$, the Fourier coefficients $c_k(\tilde \varphi)$, $k \in I_N$, are explicitly known.}

Then for all $r \in \mathbb Z$, it holds
\begin{equation*}
c_{n+r N_1}(s) = {\hat g}_n\,c_{n+r N_1}(\tilde \varphi) = \left\{ \begin{array}{ll}
c_n (f)\,\frac{c_{n+r N_1}(\tilde \varphi)}{c_n(\tilde \varphi)} &\quad n \in I_N\,,\\ [1ex]
0 & \quad n \in I_{N_1} \setminus I_N\,.
\end{array}\right.
\end{equation*}
In particular, we see that $c_n(s) = c_n(f)$ for all $n \in I_N$ and $c_n(s) = 0$ for all $n \in I_{N_1} \setminus I_N$. Substituting $k = n + r N_1$ with $n\in I_{N_1}$
and $r \in {\mathbb Z}$, we obtain
\begin{eqnarray}
s(x) - f(x) &=& \sum_{k \in \mathbb Z \setminus I_{N_1}} c_k(s)\,{\mathrm e}^{2\pi {\mathrm i}\,k x} = \sum_{n\in I_{N_1}} \sum_{r \in {\mathbb Z}\setminus \{0\}} c_{n+r N_1}(s)\,{\mathrm e}^{2\pi {\mathrm i}\,(n + r N_1)\,x}\nonumber \\
&=& \sum_{n\in I_N} c_n(f)\,\Big(\sum_{r \in {\mathbb Z}\setminus \{0\}} \frac{c_{n+r N_1}(\tilde \varphi)}{c_n(\tilde \varphi)}\,{\mathrm e}^{2\pi {\mathrm i}\,(n + r N_1)\,x} \Big)\nonumber \\
&=& \sum_{n\in I_N} c_n(f)\,\Big(\sum_{r \in {\mathbb Z}\setminus \{0\}} \frac{{\hat \varphi}(n+r N_1)}{{\hat \varphi}(n)}\,{\mathrm e}^{2\pi {\mathrm i}\,(n + r N_1)\,x} \Big)\,.\label{eq:s-fnew}
\end{eqnarray}

Let $A(\mathbb T)$ be the \emph{Wiener algebra} of all 1-periodic functions $g\in L_1(\mathbb T)$ with the property
$$
\sum_{k\in \mathbb Z} |c_k(g)| < \infty\,.
$$
Then
$$
\| g \|_{A(\mathbb T)} :=\sum_{k\in \mathbb Z} |c_k(g)|
$$
is the norm of $A(\mathbb T)$. Obviously, we have $A(\mathbb T) \subset C(\mathbb T)$, where $C(\mathbb T)$ denotes the Banach space of all 1-periodic, continuous functions with the uniform norm
$$
\| g \|_{C(\mathbb T)} :=\max_{x\in \mathbb T} |g(x)|\,.
$$
Since $x_j \in[-\frac{1}{2},\,\frac{1}{2})$, $j\in I_M$, are arbitrary nodes, we have
$$
|s(x_j) - f(x_j)| \le \| s - f\|_{C(\mathbb T)}\,.
$$
Therefore we measure the error of NFFT $\| s - f\|_{C(\mathbb T)}$ in the uniform norm. As norm of the 1-periodic trigonometric polynomial \eqref{eq:f(x)} we use the norm in the Wiener algebra $A(\mathbb T)$.
\medskip

{ 
\begin{Remark}\label{Remark:NFFTGauss}
\emph{The NFFT with a window function of the set $\Phi_{m,N_1}$ is \emph{simpler} than the NFFT with Kaiser--Bessel or Gaussian window function $\varphi$, since both window functions are supported on whole $\mathbb R$.
For such a window function $\varphi$, an additional step in the NFFT is necessary, where the 1-periodic function \eqref{eq:s(x)} is approximated by the 1-periodic well-localized function
$$
s_1(x) := \sum_{\ell \in I_{N_1}} g_{\ell}\, {\tilde \psi}\big( x - \frac{\ell}{N_1}\big)\,,
$$
where $\tilde \psi$ is the 1-periodization of the truncated window function (see \cite[pp. 378--381]{PlPoStTa18})
$$
\psi(x) := \left\{ \begin{array}{ll} \varphi(x) &\quad x \in [-m/N_1,\,m/N_1]\,,\\
0 &\quad x \in {\mathbb R} \setminus [-m/N_1,\,m/N_1]\,.
\end{array} \right.
$$
Thus the NFFT with a window function $\varphi$ supported on whole $\mathbb R$ requires also a truncated version of $\varphi$. In this case, the error of the NFFT is measured by $\|s_1 - f\|_{C(\mathbb T)}$.
In \cite[p. 393]{PlPoStTa18}, it is shown that the error of the NFFT with the Kaiser--Bessel window function can be estimated by
$$
\|s_1 - f\|_{C(\mathbb T)} \le 4 m^{3/2}\,{\mathrm e}^{- 2 \pi m \sqrt{1 - 1/\sigma}}\, \|f\|_{A(\mathbb T)}\,.
$$
We will see in Subsections 5.4 -- 5.6 that special window functions $\varphi \in \Phi_{m,N_1}$ possess a similar error behavior as the Kaiser--Bessel window function.}
\end{Remark}}
\medskip

We say that the window function $\varphi \in \Phi_{m,N_1}$ is \emph{convenient for} NFFT, if the $C(\mathbb T)$-{\emph{error constant}}
\begin{equation}
\label{eq:convenient1}
e_{\sigma}(\varphi) :=\sup_{N\in 2\mathbb N} e_{\sigma,N}(\varphi)
\end{equation}
with
\begin{equation*}
e_{\sigma,N}(\varphi) :=\max_{n\in I_N}  \big\| \sum_{r\in \mathbb Z\setminus \{0\}} \frac{ {\hat \varphi}(n+r N_1)}{{\hat \varphi}(n)}\,{\mathrm e}^{2\pi {\mathrm i}\, r N_1\,\cdot} \big\|_{C(\mathbb T)}\,, \quad N\in 2\mathbb N\,,
\end{equation*}
fulfills the condition $e_{\sigma}(\varphi) \ll 1$
for conveniently chosen truncation parameter $m \ge 2$ and oversampling factor $\sigma > 1$. Later in {  Theorem}, \ref{Lemma:esigmaNvarphibounded} we will show that {  under certain assumptions on $\varphi \in \Phi_{m,N_1}$} the value $e_{\sigma,N}(\varphi)$ is bounded for all $N \in \mathbb N$.
{  This $C(\mathbb T)$-error constant is motivated by techniques first used in \cite{St97} and later also in \cite{Ba20}. G. Steidl \cite{St97} has applied this technique for error estimates of NFFT with
B-spline and Gaussian window functions, respectively.}

\begin{Lemma}
\label{lemma:convenient2}
For each $N\in 2 \mathbb N$, the constant $e_{\sigma,N}(\varphi)$ of $\varphi \in \Phi_{m,N_1}$ can be represented in the equivalent form
\begin{equation}
\label{eq:convenient2}
e_{\sigma,N}(\varphi) = \max_{n\in I_N}  \big\| \frac{1}{N_1 c_n(\tilde \varphi)}  \sum_{\ell =0}^{N_1 -1} {\mathrm e}^{-2\pi {\mathrm i}\,n \ell/N_1}\, {\tilde \varphi}
(\cdot\, + \frac{\ell}{N_1}) - {\mathrm e}^{2\pi {\mathrm i}\,n \, \cdot} \big\|_{C(\mathbb T)}\,.
\end{equation}
\end{Lemma}

{\emph Proof.} By Lemma \ref{lemma2.1} we know that for all $x\in \mathbb T$ it holds
\begin{eqnarray*}
& &\Big(\sum_{r\in \mathbb Z\setminus \{0\}} c_{n+r N_1}(\tilde \varphi)\,{\mathrm e}^{2\pi {\mathrm i}\, r N_1\,x}\Big)\,{\mathrm e}^{2 \pi {\mathrm i}\,n x} =\,\sum_{r\in \mathbb Z\setminus \{0\}}  c_{n+r N_1}(\tilde \varphi)\,{\mathrm e}^{2\pi {\mathrm i}\,(n+ r N_1)\,x}\\
& & =\,\frac{1}{N_1}\,  \sum_{\ell =0}^{N_1 -1} {\mathrm e}^{-2\pi {\mathrm i}\,n \ell/N_1}\, {\tilde \varphi}
\big(x + \frac{\ell}{N_1}\big) - c_n(\tilde \varphi)\,{\mathrm e}^{2\pi {\mathrm i}\,n x}
\end{eqnarray*}
and hence
$$
\Big\|\sum_{r\in \mathbb Z\setminus \{0\}} \frac{ c_{n+r N_1}(\tilde \varphi)}{c_n(\tilde \varphi)}\,{\mathrm e}^{2\pi {\mathrm i}\,r N_1\,\cdot }\Big\|_{C(\mathbb T)} = \Big\|\frac{1}{N_1\,c_n(\tilde \varphi)}\,  \sum_{\ell =0}^{N_1 -1} {\mathrm e}^{-2\pi {\mathrm i}\,n \ell/N_1}\, {\tilde \varphi}\big(\cdot\, + \frac{\ell}{N_1}\big) - {\mathrm e}^{2\pi {\mathrm i}\,n \, \cdot}\Big\|_{C(\mathbb T)}\,.
$$
This completes the proof. \qedsymbol
\medskip

Thus the condition $e_{\sigma}(\varphi) \ll 1$ with \eqref{eq:convenient2} means that each exponential
\begin{equation*}
{\mathrm e}^{2\pi {\mathrm i}\,n \,\cdot}\,, \quad n \in I_N\,,
\end{equation*}
can be approximately reproduced by a linear combination of
shifted window functions ${\tilde \varphi}(\cdot\, + \frac{\ell}{N_1})$ with $\ell \in I_{N_1}$. In other words, the equispaced shifts $\tilde \varphi( \cdot + \frac{\ell}{N_1})$ with $\ell \in I_{N_1}$
are \emph{approximately exponential reproducing}.
For each node $x_j \in [-\frac{1}{2}, \frac{1}{2})$, $j \in I_M$, the linear combination
\begin{equation}
\label{eq:expapprox}
\frac{1}{N_1\,c_n(\tilde \varphi)}\,\sum_{\ell =0}^{N_1 -1} {\mathrm e}^{-2\pi {\mathrm i}\,n \ell/N_1}\, {\tilde \varphi}
(x_j + \frac{\ell}{N_1})
\end{equation}
has only few nonzero terms, since the support of $\varphi$ is very small for large $N_1$.
If we replace  $\exp (2\pi {\mathrm i} n x_j)$ for each $n\in I_N$ and $j\in I_M$ by the approximate value \eqref{eq:expapprox}, we compute approximate values of
$$
f(x_j) = \sum_{n\in I_N} c_n(f)\, {\mathrm e}^{2 \pi {\mathrm i}\,n x_j}\,, \quad j\in I_M\,,
$$
in the form
$$
f(x_j) \approx \frac{1}{N_1}\sum_{\ell \in I_{N_1}}\Big( \sum_{n\in I_N} \frac{c_n(f)}{c_n(\tilde\varphi)}\,{\mathrm e}^{-2\pi {\mathrm i}\,n \ell/N_1}\Big)\, {\tilde \varphi}
(x_j + \frac{\ell}{N_1})
$$
mainly by DFT. This is the key of the NFFT with nonequispaced spatial data and equispaced frequencies.
Special window functions $\varphi \in \Phi_{m,N_1}$ which are convenient for NFFT will be presented in Section \ref{Sec:SpWindow}.

\begin{Lemma}
\label{Lemma2.3}
Let $\sigma > 1$, $m \in \mathbb N \setminus \{1\}$, $N \in 2 \mathbb N$, and $N_1 =\sigma N \in 2 \mathbb N$ be given.
Further let $\varphi \in \Phi_{m,N_1}$. Let $f$ be a $1$-periodic trigonometric polynomial \eqref{eq:f(x)} and $s$ its approximating $1$-periodic function \eqref{eq:s(x)}.

Then the error of $\mathrm{NFFT}$ with nonequispaced spatial data and equispaced frequencies can be estimated by
\begin{equation}
\label{eq:|s-f|est}
\| s -f \|_{C(\mathbb T)} \le e_{\sigma}(\varphi)\, \|f\|_{A(\mathbb T)}\,.
\end{equation}
\end{Lemma}

\emph{Proof.} From \eqref{eq:s-fnew} it follows that
$$
s(x) - f(x) = \sum_{n\in I_N} c_n(f)\, \Big(\sum_{r \in {\mathbb Z}\setminus \{0\}} \frac{c_{n+r\sigma N}(\tilde \varphi)}{c_n(\tilde \varphi)}\,{\mathrm e}^{2\pi {\mathrm i}\,(n + r\sigma N)\,x}\Big)\,.
$$
Note that $c_n(\tilde \varphi) = {\hat \varphi}(n) > 0$ for $n \in I_N$ by assumption $\varphi \in \Phi_{m,N_1}$.
Then by H\"{o}lder's inequality we obtain that for all $x\in \mathbb T$ it holds
\begin{eqnarray*}
|s(x) - f(x)| &\le& \sum_{n\in I_N} |c_n(f)| \Big(\max_{n\in I_N}\big| \sum_{r \in {\mathbb Z}\setminus \{0\}} \frac{c_{n+r N_1}(\tilde \varphi)}{c_n(\tilde \varphi)}\,{\mathrm e}^{2\pi {\mathrm i}\,r N_1\,x}\big|\Big)\\
&=& \sum_{n\in I_N} |c_n(f)| \Big(\max_{n\in I_N}\big| \sum_{r \in {\mathbb Z}\setminus \{0\}} \frac{{\hat \varphi}(n+r N_1)}{{\hat \varphi}(n)}\,{\mathrm e}^{2\pi {\mathrm i}\,r N_1\,x}\big|\Big)\\
&\le&  e_{\sigma}(\varphi)\, \|f\|_{A(\mathbb T)} < \infty\,.
\end{eqnarray*}
Hence we get \eqref{eq:|s-f|est}. \qedsymbol
\medskip

\begin{Remark}
\emph{Let $\lambda \ge 0$ be fixed. We introduce the $1$-\emph{periodic Sobolev space} $H^{\lambda}(\mathbb T)$ of all $1$-periodic functions $f:\mathbb T \to \mathbb C$ which are integrable on $[0,\,1]$
and for which
$$
\|f\|_{H^{\lambda}(\mathbb T)} := \Big(\sum_{k\in \mathbb Z} |k|^{2 \lambda} |c_k(f)|^2\big)^{1/2} < \infty\,,
$$
where we {  declare} $|0|:=1$. Then $H^{\lambda}(\mathbb T)$ is a Hilbert space with the inner product
$$
\langle f,\,g \rangle_{H^{\lambda}(\mathbb T)}:= \sum_{k\in \mathbb Z} |k|^{2 \lambda} c_k(f)\,\overline{c_k(g)}\,.
$$
For $\lambda =0$, we have $H^{0}(\mathbb T)=L^{2}(\mathbb T)$. Then the Sobolev embedding theorem (see \cite[p.~142]{SaVa}) says that for $\lambda > \frac{1}{2}$ it holds
$
H^{\lambda}(\mathbb T) \subset A(\mathbb T) \subset C(\mathbb T)\,.
$
Let $f$ be a $1$-periodic, trigonometric polynomial of the form \eqref{eq:f(x)}. Then we have
\begin{eqnarray*}
\|f\|_{C(\mathbb T)} &\le& \sum_{k\in I_N} |c_k(f)| = \|f\|_{A(\mathbb T)}
= \sum_{k\in I_N} \big(|c_k(f)|\,|k|^{\lambda}\big)\,|k|^{-\lambda}\,.
\end{eqnarray*}
Applying the Cauchy--Schwarz inequality, we obtain for $\lambda > \frac{1}{2}$ that
\begin{eqnarray*}
\|f\|_{C(\mathbb T)} &\le& \|f\|_{A(\mathbb T)} \le \Big(\sum_{k\in I_N} \big(|c_k(f)|^2\,|k|^{2\lambda}\Big)^{1/2}\Big(\sum_{k\in I_N} |k|^{-2\lambda}\Big)^{1/2}\\
&\le& \|f\|_{H^{\lambda}(\mathbb T)} \Big(1 + 2 \sum_{k=1}^{\infty} \frac{1}{k^{2\lambda}}\Big)^{1/2}\,.
\end{eqnarray*}
Using the \emph{Riemann zeta function}
$
\zeta(2\lambda) := \sum_{k=1}^{\infty} \frac{1}{k^{2\lambda}}$, $\lambda > \frac{1}{2},
$
we obtain the following inequality
$$
\|f\|_{C(\mathbb T)} \le \|f\|_{A(\mathbb T)} \le \|f\|_{H^{\lambda}(\mathbb T)} \sqrt{1 + 2 \zeta(2\lambda)}\,.
$$
Thus under the assumptions of Lemma \ref{Lemma2.3}, the error of NFFT with nonequispaced spatial data and equispaced frequencies can be estimated by
$$
\| s - f \|_{C(\mathbb T)} \le e_{\sigma}(\varphi)\,\sqrt{1 + 2 \zeta(2\lambda)}\, \|f\|_{H^{\lambda}(\mathbb T)}\,.
$$
Note that}
$$
\zeta(2) = \frac{\pi^2}{6},\;\zeta(4) = \frac{\pi^4}{90},\;\zeta(6) = \frac{\pi^6}{945},\;\zeta(8) = \frac{\pi^8}{9450},\; \zeta(10) = \frac{\pi^{10}}{93555} < 1.000995\,.
$$
\end{Remark}

\subsection{NFFT with nonequispaced frequencies and equispaced spatial data}

The NFFT \emph{with nonequispaced frequencies and equispaced spatial data} or \emph{transposed} NFFT evaluates the exponential sums
\begin{equation}
\label{eq:sk:=}
s_k := \sum_{j\in I_M} f_j\,{\mathrm e}^{2 \pi {\mathrm i}\,k x_j}\,, \quad k\in I_N\,,
\end{equation}
for arbitrary given coefficients $f_j \in \mathbb C$ and {  nonequispaced frequencies} $x_j\in \big[-\frac{1}{2},\,\frac{1}{2}\big)$, $j\in I_M$. Assume that the window function $\varphi\in \Phi_{m,N_1}$ is convenient for NFFT.
Introducing the 1-periodic function
$$
g(x) := \sum_{j\in I_M} f_j\,{\tilde \varphi}(x_j + x)\,,
$$
the Fourier coefficients of $g$ read as follows
$$
c_k(g) = \int_0^1 g(t)\,{\mathrm e}^{- 2\pi {\mathrm i}\,k t}\,{\mathrm d}t = \Big(\sum_{j \in I_M} f_j\,{\mathrm e}^{2 \pi {\mathrm i}\,k x_j}\Big)\,c_k(\tilde \varphi) = s_k\,c_k(\tilde \varphi)\,, \quad k \in \mathbb Z\,.
$$
Using the trapezoidal rule, we {  approximate} $c_k(g)$ by
$$
{\widetilde c}_k(g) := \frac{1}{N_1}\, \sum_{\ell\in I_{N_1}} \sum_{j\in I_M} f_j\,{\tilde \varphi}\big(x_j + \frac{\ell}{N_1}\big)\,{\mathrm e}^{- 2\pi {\mathrm i}\,k \ell/N_1}\,.
$$
{  Note that ${\widetilde c}_k(g)$, $k \in I_N$, can be efficiently computed by FFT, for details see \cite[p. 382]{PlPoStTa18}.}
Then the results of this NFFT with nonequispaced frequencies and equispaced spatial data are the values
\begin{equation}
\label{eq:hatsk:=}
{\widetilde s}_k := \frac{{\widetilde c}_k(g)}{c_k(\tilde \varphi)}\,, \quad k \in I_N\,.
\end{equation}
It is interesting that the same $C(\mathbb T)$-error constant \eqref{eq:convenient1} appears in an error estimate of the NFFT with nonequispaced frequencies and equispaced spatial data too.

\begin{Lemma}\label{lemma:convenientNFFT}
Let $\sigma > 1$, $m \in \mathbb N \setminus \{1\}$, $N,\,M \in 2 \mathbb N$, and $N_1 = \sigma N \in 2 \mathbb N$. Further let $\varphi \in \Phi_{m,N_1}$. For given $f_j \in \mathbb C$ and nonequispaced
frequencies $x_j \in \big[- \frac{1}{2},\,
\frac{1}{2}\big)$, $j \in I_M$, we consider the exponential sums \eqref{eq:sk:=} and the related approximations \eqref{eq:hatsk:=}.

Then the error of $\mathrm{NFFT}$ with nonequispaced frequencies and equispaced spatial data can be estimated by
$$
\max_{k\in I_N} | s_k - {\widetilde s}_k| \le e_{\sigma}(\varphi)\, \sum_{j\in I_M} |f_j|\,.
$$
\end{Lemma}

\emph{Proof.} For each $k \in I_N$ we have
\begin{eqnarray*}
|s_k - {\widetilde s}_k| &=& \frac{1}{c_k(\tilde \varphi)}\,|c_k(g) - {\widetilde c}_k(g)|\\ [1ex]
&=& \Big|\sum_{j\in I_M} f_j\,\Big( {\mathrm e}^{2\pi {\mathrm i}\,k x_j} - \frac{1}{N_1 c_k(\tilde \varphi)}\,\sum_{\ell\in I_{N_1}} {\tilde \varphi}\big(x_j + \frac{\ell}{N_1}\big)\,{\mathrm e}^{- 2\pi {\mathrm i}\,k \ell/N_1}\Big)\Big|\,.
\end{eqnarray*}
From H\"{o}lder's inequality and Lemma \ref{lemma:convenient2} it follows that
\begin{eqnarray*}
|s_k - {\widetilde s}_k| &\le& \max_{j\in I_M} \big| {\mathrm e}^{2\pi {\mathrm i}\,k x_j} - \frac{1}{N_1 c_k(\tilde \varphi)}\,\sum_{\ell\in I_{N_1}} {\tilde \varphi}\big(x_j + \frac{\ell}{N_1}\big)\,{\mathrm e}^{- 2\pi {\mathrm i}\,k \ell/N_1}\big|\,\sum_{j\in I_M} |f_j|\\
&\le& e_{\sigma, N}(\varphi)\,\sum_{j\in I_M} |f_j| \le e_{\sigma}(\varphi)\,\sum_{j\in I_M} |f_j|\,.
\end{eqnarray*}
This completes the proof. \qedsymbol

\section{Auxiliary estimates}\label{Sec:FT}

In our study we use later the following

\begin{Lemma}
\label{Lemma3.7}
For $-1 < u < 1$ and $\mu > 1$ it holds
$$
\sum_{r \in {\mathbb Z}\setminus \{0,\,\pm 1\}} |u + r|^{-\mu} \le \frac{2}{\mu -1}\, (1 - |u|)^{1-\mu}\,.
$$
\end{Lemma}

\emph{Proof.} For $-1 < u < 1$, $r\in \mathbb N$, and $\mu > 1$ we have
\begin{equation}
\label{eq:|u+r|}
| u + (-1)^j\,r|^{-\mu} \le (r - |u|)^{-\mu}\,, \quad j= 0,\,1\,.
\end{equation}
Using \eqref{eq:|u+r|}, the series can be estimated as follows
$$
\sum_{r=2}^{\infty} |u + (-1)^j\,r|^{-\mu} \le  \sum_{r=2}^{\infty} (r - |u|)^{-\mu}\,, \quad j= 0,\,1\,.
$$
Hence it follows by the integral test for convergence that
\begin{eqnarray*}
\sum_{r \in {\mathbb Z}\setminus \{0, \pm 1\}} |u + r|^{-\mu} &=& \sum_{r=2}^{\infty} |u +r|^{-\mu} + \sum_{r=2}^{\infty} |u -r|^{-\mu}
\le 2\,\sum_{r=2}^{\infty} (r - |u|)^{-\mu}\\ [1ex]
&\le&  2\, \int_1^{\infty}(x - |u|)^{-\mu}\, {\mathrm d}x = \frac{2}{\mu -1}\,(1 -|u|)^{1-\mu}\,. \quad \qedsymbol
\end{eqnarray*}
\medskip

For fixed $\mu \ge 0$, the $\mu${\emph{th Bessel function of first kind}} is defined by
\begin{equation*}
J_{\mu}(x) := \sum_{k=0}^{\infty}\frac{(-1)^k}{2^{\mu + 2k}\, k!\, \Gamma(\mu+k+1)}\, x^{\mu+2k}\,,\quad x \ge 0\,,
\end{equation*}
so that in particular
$$
J_m(x) := \sum_{k=0}^{\infty}\frac{(-1)^k}{2^{m + 2k}\, k!\, (m+k)!}\, x^{m+2k}\,, \quad m=0,\,1,\,\ldots\,.
$$
For fixed $\mu \ge 0$, the $\mu${\emph{th modified Bessel function of first kind}} is defined by
\begin{equation*}
I_\mu(x) := \sum_{k=0}^{\infty}\frac{1}{2^{\mu + 2k}\, k!\, \Gamma(\mu+k+ 1)}\, x^{\mu+2k}\,,\quad x \ge 0\,,
\end{equation*}
so that in particular
\begin{equation*} 
I_m(x) := \sum_{k=0}^{\infty}\frac{1}{2^{m + 2k}\, k!\, (m+k)!}\, x^{m+2k}\,, \quad m=0,\,1,\,\ldots\,.
\end{equation*}
For the properties of Bessel functions we refer to \cite[pp.~355--478]{abst} and \cite{Watson}. In particular, these Bessel functions possess the following asymptotic behaviors for $x\to \infty$ (see \cite[pp.~364, 377]{abst}),
\begin{eqnarray}
J_{\mu}(x) &\sim& \sqrt{\frac{2}{\pi x}}\, \cos\big(x - \frac{\mu \pi}{2} - \frac{\pi}{4}\big)\,\big(1 + {\mathcal O}(x^{-1})\big)\,, \label{eq:Jmusim}\\
I_{\mu}(x) &\sim& \frac{1}{\sqrt{2\pi x}}\, {\mathrm e}^x\,\big(1  - \frac{4 \mu^2 - 1}{8 x} + {\mathcal O}(x^{-2})\big)\,. \label{eq:Imusim}
\end{eqnarray}
Here we are interested in \emph{explicit error estimates} for NFFT with compactly supported, continuous window function. For this purpose, we need explicit bounds for the Bessel functions instead of
the asymptotic formulas \eqref{eq:Jmusim} and \eqref{eq:Imusim}.

\begin{Lemma}
For fixed $\mu > \frac{1}{2}$ and all $x\ge 0$, it holds
\begin{equation}
\label{eq:|Jmu(x)|}
\big| x^2 - \mu^2 + \frac{1}{4}\big|^{1/4} \,|J_{\mu}(x)| < \sqrt{\frac{2}{\pi}}\,,
\end{equation}
where $\sqrt{\frac{2}{\pi}}$ is the best possible upper bound.
In particular for $\mu =1$ and $x \ge 6$, we have
\begin{equation*}
|J_1(x)| < \frac{1}{\sqrt x}\,.
\end{equation*}
For $\mu =\frac{5}{2}$ and $x \ge 6$, we have
\begin{equation}
\label{eq:|J5/2(x)|}
|J_{5/2}(x)| < \frac{1}{\sqrt x}\,.
\end{equation}
For $\mu =3m$ with $m \in \mathbb N \setminus \{1\}$ and $x \ge \pi m$, we have
\begin{equation}
\label{eq:|J3m(x)|}
|J_{3m}(x)| < \frac{3}{2\,\sqrt x}\,.
\end{equation}
\end{Lemma}

\emph{Proof.} The inequality \eqref{eq:|Jmu(x)|} was shown in \cite{Kra}.
For $\mu =1$, the inequality \eqref{eq:|Jmu(x)|} means that
$$
\big| x^2 -  \frac{3}{4}\big|^{1/4} \,|J_{1}(x)| < \sqrt{\frac{2}{\pi}}\,.
$$
Thus for all $x \ge 6$ we have
$$
c\,\sqrt x \le \big| x^2 -  \frac{3}{4}\big|^{1/4}
$$
with the constant
$$
c = \min_{x \ge 6} \frac{ \big| x^2 -  \frac{3}{4}\big|^{1/4}}{\sqrt x} = \min_{x \ge 6} \big(1- \frac{3}{4x^2}\big)^{1/4}= \big(1 - \frac{1}{48}\big)^{1/4} = 0.994750\ldots
$$
such that
$$
|J_1(x)| \le \frac{1}{0.994750}\,\sqrt{\frac{2}{\pi x}} < \frac{1}{\sqrt x}\,.
$$
Similarly, one can show the inequalities \eqref{eq:|J5/2(x)|} and \eqref{eq:|J3m(x)|}. \qedsymbol

\begin{Lemma}
Let $\mu \ge \frac{1}{2}$ and $x_0 >0$ be given. Then for all $x\ge x_0$, it holds
\begin{equation}
\label{eq:Imu(x)estim}
\sqrt{2 \pi x_0}\,{\mathrm e}^{-x_0}\,I_{\mu}(x_0) \le \sqrt{2 \pi x}\,{\mathrm e}^{-x}\,I_{\mu}(x) < 1\,.
\end{equation}
\end{Lemma}

\emph{Proof.} By { \cite[Inequality (2.6)]{Bar}} one knows that for fixed $\mu \ge \frac{1}{2}$ and arbitrary $x$, $y\in (0,\,\infty)$ with $x<y$ it holds
$$
\frac{I_{\mu}(x)}{I_{\mu}(y)} < \sqrt{\frac{y}{x}}\,{\mathrm e}^{x-y}\,,
$$
i.e.,
$$
\sqrt{2\pi x}\,{\mathrm e}^{-x}\,I_{\mu}(x) < \sqrt{2\pi y}\,{\mathrm e}^{-y}\,I_{\mu}(y)\,.
$$
Hence the function $f(x):= \sqrt{2\pi x}\,{\mathrm e}^{-x}\,I_{\mu}(x)$ is strictly increasing on $(0,\,\infty)$. Further by \eqref{eq:Imusim} we have
$$
\lim_{y \to \infty} \sqrt{2\pi y}\,{\mathrm e}^{-y}\,I_{\mu}(y) = 1\,.
$$
This implies the inequality \eqref{eq:Imu(x)estim}. \qedsymbol

\section{Paley--Wiener theorem in a Sobolev space}\label{Sec:Paley-Wiener}

The main tool of this approach is the study of the Fourier transform ${\hat \varphi}(v)$ of $\varphi \in \Phi_{m,N_1}$ for $|v| \to \infty$. A rapid decay of $\hat \varphi$ is essential for a small
$C(\mathbb T)$-error constant \eqref{eq:convenient1}. From Fourier analytical point of view, it is very interesting to discuss the relation to the
known Theorem of Paley--Wiener (see \cite[pp.~12--13]{PaWi}), since this result characterizes the behavior of Fourier transforms of functions $\varphi$ which vanish outside the open
interval $I =\big(-\frac{m}{N_1},\,\frac{m}{N_1}\big)$. The smoothness of the restricted window function $\varphi|_I$ determines the decay of ${\hat \varphi}(v)$ for $|v| \to \infty$.

For simplicity,
we denote functions defined on $I$ by $\varphi$ too.
By $H^k(I)$, $k\in {\mathbb N}$, we denote the Sobolev space of all functions $\varphi \in L^2(I)$
with $D^j \varphi \in L^2(I)$, $j = 1,\ldots, k$, where $D^j \varphi$ is the $j$th weak derivative of $\varphi$. Then $H^k(I)$ is a Hilbert space with the Sobolev norm
$$
\|\varphi\|_{H^k(I)} := \Big(\sum_{j=0}^k \| D^j \varphi\|_{L^2(I)}^2\Big)^{1/2}\,.
$$
Further we define the Sobolev space $H_0^k(I)$ as the closure in $H^k(I)$ of the space $C_0^{\infty}(I)$ which consists of all infinitely differentiable, locally supported functions.
{  Then $H_0^k(I)$ is the set of all functions $\varphi \in H^k(I)$ with vanishing one-sided derivatives $\varphi^{(j)}\big(\pm \frac{m}{N_1}\big) = 0$, $j =0,\ldots k-1$.}
We present the following modification of the Theorem of Paley--Wiener (cf. \cite[pp.~12--13]{PaWi} or \cite[Vol.~II, pp.~272--274]{Zy}):
\medskip

\begin{Theorem}
\label{Theorem:Paley-Wiener}
For given $\varphi \in H_0^k(I)$ with $k\in {\mathbb N}$, the function
$$
f(z) := \int_I \varphi(t)\,{\mathrm e}^{- 2 \pi{\mathrm i}z t}\,{\mathrm d}t\,, \quad z= x + {\mathrm i}\,y \in \mathbb C\,,
$$
is entire and has the following properties:\\
$1$. For real variable $x$, it holds $f(x) \in L^2(\mathbb R)$ and $x^k\,f(x) \in L^2(\mathbb R)$.\\
$2$. For all $z \in \mathbb C$, there exist a positive constant $\gamma_k(\varphi)$ such that
$$
|f(z)| \le \gamma_k(\varphi)\, (1 + |2\pi z|)^{-k}\, {\mathrm e}^{2\pi m\,|{\mathrm{Im}}z|/N_1}\,.
$$
Conversely, if an entire function $f$ satisfies the conditions $1$. and $2$., then the function
$$
\varphi(x) :=  \int_{\mathbb R} f(t)\, {\mathrm e}^{2 \pi {\mathrm i}x t}\,{\mathrm d}t\,, \quad x \in \mathbb R\,,
$$
has the properties $\varphi |_ I \in H^k(I)$ and $\varphi |_{ {\mathbb R} \setminus I} = 0$. Note that
$$
\int_{\mathbb R} f(t)\, {\mathrm e}^{2 \pi{\mathrm i}x t}\,{\mathrm d}t :=  \lim_{R\to \infty} \int_{-R}^R f(t)\, {\mathrm e}^{2 \pi{\mathrm i}x t}\,{\mathrm d}t
$$
means the limit in $L^2(\mathbb R)$ and that $\varphi |_{ {\mathbb R} \setminus I} = 0$ means $\varphi(x) = 0$ for almost all $x\in {\mathbb R} \setminus I$.
\end{Theorem}
\medskip

{\em Proof.} Let $\varphi \in H_0^k(I)$ be given, i.e., $D^j\varphi \in L^2(I)$, $j =0,\ldots, k$, and
\begin{equation}
\label{eq:derivativesphi}
\varphi^{(j)}\big(\pm \frac{m}{N_1}\big) = 0\,, \quad j = 0,\, \ldots ,\,k-1\,.
\end{equation}
For arbitrary $z= x + {\mathrm i}\,y \in
\mathbb C$ and $t \in \bar I$ we have
$$
|{\mathrm e}^{- 2 \pi{\mathrm i}z t}| = |{\mathrm e}^{- 2 \pi{\mathrm i}x t}|\,{\mathrm e}^{2 \pi yt}\le {\mathrm e}^{2 \pi |y|\,|t|} \le {\mathrm e}^{2 \pi m\,|{\mathrm{Im}}z|/N_1}
$$
and hence
\begin{equation}\label{1}
|f(z)| = \big| \int_I \varphi(t)\,{\mathrm e}^{- 2 \pi{\mathrm i}z t}\,{\mathrm d}t\big| \le \big(\int_I |\varphi(t)|\,{\mathrm d}t\big)\,{\mathrm e}^{2 \pi m\,|{\mathrm{Im}}z|/N_1}\,.
\end{equation}
Since by the Schwarz inequality
\begin{equation} \label{1a}
\|\varphi\|_{L^1(I)} = \int_I |\varphi(t)|\cdot 1\,{\mathrm d}t \le \|\varphi\|_{L^2(I)}\,\|1\|_{L^2(I)} = \sqrt{\frac{2m}{N_1}}\,\|\varphi\|_{L^2(I)} < \infty\,,
\end{equation}
we have $\varphi \in L^1(I)$ too. Since $D^j\varphi \in L^2(I)$ for $j = 1,\,\ldots, \,k$, we have $D^j\varphi \in L^1(I)$ analogously to \eqref{1a}.

Obviously, $f$ is an entire function, because for each $z\in \mathbb C$ and $j \in \mathbb N$ it holds
$$
f^{(j)}(z) = \int_I \varphi(t)\,(- 2 \pi{\mathrm i}\,t)^j\,{\mathrm e}^{- 2 \pi{\mathrm i}z t}\,{\mathrm d}t = (- 2 \pi {\mathrm i})^j \,\int_I t^j\,\varphi(t)\,{\mathrm e}^{- 2 \pi{\mathrm i}z t}\,{\mathrm d}t\,.
$$
Using \eqref{eq:derivativesphi}, repeated integration by parts applied to the function $\varphi$ gives for $j=1,\ldots, k$ the equalities
$$
\int_I \varphi(t)\,\Big(\frac{{\mathrm d}^j}{{\mathrm d}t^j}\,{\mathrm e}^{- 2 \pi{\mathrm i}z t}\Big)\,{\mathrm d}t = (-1)^j\, \int_I (D^j \varphi)(t)\,{\mathrm e}^{- 2 \pi{\mathrm i}z t}\,{\mathrm d}t
$$
such that
\begin{equation}
\label{2}
(- 2 \pi{\mathrm i}z)^j\,f(z) =
\int_I \varphi(t)\,\Big(\frac{{\mathrm d}^j}{{\mathrm d}t^j}\,{\mathrm e}^{- 2 \pi{\mathrm i}z t}\Big)\,{\mathrm d}t = (-1)^j\,\int_I (D^j \varphi)(t)\,{\mathrm e}^{- 2 \pi{\mathrm i}z t}\,{\mathrm d}t\,.
\end{equation}
Hence we obtain
\begin{equation}
\label{3}
|2 \pi z|^j \, |f(z)| \le \big(\int_I |(D^j \varphi)(t)|\,{\mathrm d}t\big)\, {\mathrm e}^{2 \pi m\,|{\mathrm{Im}}z|/N_1}\,.
\end{equation}
From \eqref{1} and \eqref{3} it follows that for all $z \in \mathbb C$ it holds
$$
(1 + |2\pi z|)^k\, |f(z)| \le \gamma_k(\varphi)\, {\mathrm e}^{2\pi m\,|{\mathrm{Im}}z|/N_1}
$$
with the positive constant
$$
\gamma_k(\varphi):= {k \choose 0}\,\int_I |\varphi(t)|\,{\mathrm d}t + {k \choose 1}\,\int_I |(D\varphi)(t)|\,{\mathrm d}t + \ldots + {k \choose k}\,\int_I |(D^k \varphi)(t)|\,{\mathrm d}t < \infty\,.
$$
By \eqref{2} with $j=k$, we have for $x \in \mathbb R$ the equalities
$$
(- 2 \pi{\mathrm i}x)^k\,f(x)  = (-1)^k\,\int_I (D^k \varphi)(t)\,{\mathrm e}^{- 2 \pi{\mathrm i}x t}\,{\mathrm d}t = (-1)^k\,\int_{\mathbb R} (D^k \varphi)(t)\,{\mathrm e}^{- 2 \pi{\mathrm i}x t}\,{\mathrm d}t
$$
such that $x^k\,f(x) \in L^2(I)$ by the Theorem of Plancherel.
\medskip

Now we assume that an entire function $f$ with the properties 1. and 2. is given. Then especially we have $f|_{\, \mathbb R} \in L^2(\mathbb R)$ and
$$
|f(z)| \le \gamma_k(\varphi) \,  {\mathrm e}^{2\pi m\,|{\mathrm{Im}}z|/N_1}\,, \quad z \in \mathbb C\,.
$$
By the original Theorem of Paley--Wiener, the function
$$
\varphi(x) := \int_{\mathbb R} f(t)\,{\mathrm e}^{2 \pi{\mathrm i}t x}\, {\mathrm d}t \in L^2(\mathbb R)
$$
vanishes for almost all $x \in {\mathbb R}\setminus I$ such that $\varphi |_I \in L^2(I)$. For $t\in \mathbb R$, from $f(t) \in L^2(\mathbb R)$ and $t^k\,f(t) \in L^2(\mathbb R)$ it follows $t^j\,f(t) \in L^2(\mathbb R)$, $j=1,\ldots,k-1$,
since $|t^j\,f(t)|^2 \le |f(t)|^2$ for $t\in [-1,\,1]$ and  $|t^j\,f(t)|^2 \le |t^k\,f(t)|^2$ for $t\in \mathbb R \setminus [-1,\,1]$.

Each function in $L^2(\mathbb R)$ generates a tempered distribution. By the differentiation property of the Fourier transform of tempered distributions we conclude that
the $j$th weak derivative of $\varphi$ exists almost everywhere and that
$$
(D^j \varphi)(x) = \int_{\mathbb R} (2 \pi {\mathrm i}t)^j\,f(t)\,{\mathrm e}^{2 \pi{\mathrm i}t x}\, {\mathrm d}t \in L^2(\mathbb R)\,, \quad j = 1,\ldots , k.
$$
Now we have to show that $D^j \varphi |_{{\mathbb R}\setminus I} = 0$.
From property 2. of the entire function $f$ it follows that the entire function $(2\pi {\mathrm i}z)^j\,f(z)$ with $j = 1,\ldots, k$ fulfills the inequality
$$
\big|(2 \pi {\mathrm i}z)^j\,f(z)\big| = |2 \pi z|^j\,|f(z)| \le \gamma_k(\varphi)\,\frac{|2\pi z|^j}{(1+|2 \pi z|)^k}\, {\mathrm e}^{2 \pi m\,|{\mathrm{Im}}z|/N_1} \le \gamma_k(\varphi)\,{\mathrm e}^{2 \pi m\,|{\mathrm{Im}}z|/N_1}
$$
for all $z\in \mathbb C$.
Applying the original Theorem of Paley--Wiener to the entire function $(2 \pi{\mathrm i}z)^j\,f(z)$, we see that
$$
\psi_j(x) := \int_{\mathbb R}(2 \pi {\mathrm i}t)^j\,f(t)\,{\mathrm e}^{2 \pi{\mathrm i}t x}\, {\mathrm d}t
$$
vanishes almost everywhere in $\mathbb R \setminus I$. Since $\psi_j(x) = (D^j \varphi)(x) =0$ almost everywhere, we have $(D^j \varphi)(x) =0$ for almost all $x \in \mathbb R \setminus I$.
Hence we have $D^j \varphi |_{I}  \in L^2(I)$ for $j=1,\ldots,k$ such that $\varphi |_I \in H^k(I)$. This completes the proof. \qedsymbol
\medskip

{ 
Unfortunately, the assumption $\varphi \in H_0^k(I)$ is too strong for the most popular window function.

\begin{Example}
\label{Example:Paley-Wiener}
The triangular function $\varphi(x) := 1 - \frac{N_1}{m}\,|x|$, $x \in I$, belongs to the Sobolev space $H_0^1(I)$, but doesn't belong to $H_0^2(I)$, since it holds
$\varphi\big(\pm \frac{m}{N_1}\big) = 0$ and $\varphi'\big(\pm \frac{m}{N_1}\big) \neq 0$. In this case, we obtain that for $z\in \mathbb C$,
\begin{eqnarray*}
f(z) &=& \int_{-m/N_1}^{m/N_1} \varphi(t)\,{\mathrm e}^{-2 \pi {\mathrm i}\,z t}\,{\mathrm d}t = 2\,\int_0^{m/N_1} \big(1 - \frac{N_1}{m}\,t\big)\, \cos(2\pi z t)\,{\mathrm d}t\\
&=& \frac{2m}{N_1}\,\int_0^1 (1 - s)\,\cos \big(\frac{2 \pi m}{N_1}\, z s \big)\, {\mathrm d}s = \frac{m}{N_1}\,\big(\mathrm{sinc}\,\frac{\pi m z}{N_1}\big)^2\,.
\end{eqnarray*}
Thus Theorem \ref{Theorem:Paley-Wiener} results in $k = 1$. Otherwise we observe a quadratic decay of the Fourier transform $\hat \varphi$, since for $x \in \mathbb R \setminus \{0\}$,
$$
|f(x)| = |{\hat \varphi}(x)| \le \frac{N_1}{m \pi^2}\,|x|^{-2}\,.
$$
\end{Example}

Hence we will present a better method in the following Theorem \ref{Lemma:esigmaNvarphibounded}.
}

\section{Special  window functions}\label{Sec:SpWindow}

In this section, we determine upper bounds of the error constant \eqref{eq:convenient1} for various special window functions $\varphi \in \Phi_{m,N_1}$ by two methods.
If the series
\begin{equation*}
\sum_{n\in \mathbb Z} |{\hat \varphi}(n)| < \infty\,,
\end{equation*}
then by \eqref{eq:convenient1} we have that
\begin{equation}
\label{eq:boundesimple}
e_{\sigma}(\varphi) \le \sup_{N\in 2\mathbb N} \Big(\max_{n\in I_N} \, \sum_{r \in {\mathbb Z}\setminus \{0\}} \frac{|{\hat \varphi}(n+r N_1)|}{{\hat \varphi}(n)}\Big)\,.
\end{equation}
This technique can be applied for (modified) B-spline window functions.\\

For the algebraic, Bessel, $\sinh$-type, and modified $\cosh$-type window functions, we use the following argument. By H\"older's inequality it follows
from \eqref{eq:convenient1} that
$$
e_{\sigma}(\varphi) \le \sup_{N\in 2\mathbb N}\big(\min_{n\in I_N} {\hat\varphi}(n)\big)^{-1}\,\max_{n\in I_N}  \big\|\sum_{r \in {\mathbb Z}\setminus \{0\}} {\hat \varphi}(n + r N_1)\,{\mathrm e}^{2 \pi {\mathrm i}\,N_1 r \, \cdot}\big\|_{C(\mathbb T)}\,.
$$
Thus we show that the minimum of all ${\hat \varphi}(n)$ for relatively low frequencies $n\in I_N$ is equal to ${\hat\varphi}(N/2)$ and
that the series
$$
\sum_{r \in {\mathbb Z}\setminus \{0\}} |{\hat \varphi}(n+r N_1)|
$$
is bounded for each $n\in I_N$. For this we have to estimate the Fourier transform ${\hat \varphi}(n)$ for sufficiently large frequencies
$|n| \ge N_1 -\frac{N}{2}$ very carefully. Thus  this results in
\begin{equation*}
e_{\sigma}(\varphi) \le \sup_{N\in 2\mathbb N}\Big(\frac{1}{{\hat\varphi}(N/2)}\,\max_{n\in I_N} \sum_{r \in \mathbb Z \setminus \{0\}} |{\hat \varphi}(n + r N_1)|\Big)\,.
\end{equation*}

Now we use the special structure of the given window function $\varphi$.
Let $m \in {\mathbb N}\setminus \{1\}$ and $\sigma > 1$ be given.
Assume that an even, continuous function $\varphi_0: {\mathbb R} \to [0,\,1]$ with $\mathrm{supp}\,\varphi_0 =
[-1,\,1]$ has following properties: $\varphi_0(0) = 1$, $\varphi_0|_{[0,\,1]}$ is decreasing, and its
restricted Fourier transform ${\hat \varphi}_0|_{[0, \,m/(2\sigma)]}$ is positive and decreasing.
Let $N \in 2{\mathbb N}$ with $N_1 = \sigma N \in 2{\mathbb N}$ be given. Then the scaled function
\begin{equation}
\label{eq:scaledvarphi0}
\varphi(x) := \varphi_0\big(\frac{N_1 x}{m}\big)\,, \quad x \in \mathbb R\,,
\end{equation}
is a continuous window function of the set $\Phi_{m,N_1}$.

\begin{Theorem}
\label{Lemma:esigmaNvarphibounded}
Let $\sigma > 1$, $m \in \mathbb N \setminus \{1\}$, $N \in 2 \mathbb N$, and $N_1 = \sigma N \in 2 \mathbb N$ be given. Further let $\varphi \in \Phi_{m,N_1}$ be a scaled version \eqref{eq:scaledvarphi0} of $\varphi_0$. Assume that
the Fourier transform ${\hat \varphi}_0$ fulfills the decay condition
$$
|{\hat \varphi}_0(v)| \le \left\{ \begin{array}{ll} c_1 & \quad |v| \in \big[m \big(1 - \frac{1}{2\sigma}\big),\,m \big(1 + \frac{1}{2\sigma}\big)\big]\,,\\ [1ex]
 c_2 \,|v|^{-\mu} & \quad |v| \ge m \big(1+ \frac{1}{2\sigma}\big)\,,
\end{array} \right.
$$
with certain constants $c_1 > 0$, $c_2>0$, and $\mu > 1$.

Then the constant $e_{\sigma,N}(\varphi)$ is bounded for all $N \in 2{\mathbb N}$ by
$$
e_{\sigma,N}(\varphi) \le \frac{1}{{\hat \varphi}_0\big(\frac{m}{2 \sigma}\big)} \Big[2 c_1 + \frac{2 c_2}{(\mu -1)\,m^{\mu}}\big(1 - \frac{1}{2\sigma}\big)^{1-\mu}\Big]\,.
$$
Further, the $C(\mathbb T)$-error constant $e_{\sigma}(\varphi)$ of the window function \eqref{eq:scaledvarphi0} has the upper bound
\begin{equation}
\label{eq:boundesigmascaledvarphi}
e_{\sigma}(\varphi) \le \frac{1}{{\hat \varphi}_0\big(\frac{m}{2 \sigma}\big)} \Big[2 c_1 + \frac{2 c_2}{(\mu -1)\,m^{\mu}}\big(1 - \frac{1}{2\sigma}\big)^{1-\mu}\Big]\,.
\end{equation}
\end{Theorem}

\emph{Proof.} Note that it holds
$$
|{\hat \varphi}_0(v)| \le \int_{\mathbb R} \varphi_0(x)\,{\mathrm d}x = {\hat \varphi}_0(0)\,, \quad v\in \mathbb R\,.
$$
By the scaling property of the Fourier transform, we have
$$
{\hat \varphi}(v) = \int_{\mathbb R} \varphi(x)\,{\mathrm e}^{-2 \pi {\mathrm i}v x}\,{\mathrm d}x = \frac{m}{N_1}\,{\hat \varphi}_0\big(\frac{m v}{N_1}\big)\,, \quad v \in \mathbb R\,.
$$
For all $n \in I_N$ and $r\in {\mathbb Z}\setminus \{0,\,\pm 1\}$, we obtain
$$
\big| \frac{m n}{N_1} + m r \big| \ge m \big(2 - \frac{1}{2 \sigma}\big) > m \big(1 + \frac{1}{2 \sigma}\big)
$$
and hence
$$
|{\hat \varphi}(n + r N_1)| = \frac{m}{N_1}\,\big|{\hat \varphi}_0\big(\frac{m n}{N_1} + m r\big)\big| \le \frac{m\,c_2}{m^{\mu} N_1}\,\big|\frac{n}{N_1} + r\big|^{-\mu}\,.
$$
From Lemma \ref{Lemma3.7} it follows that for fixed $u =\frac{n}{N_1}\in \big[- \frac{1}{2\sigma},\, \frac{1}{2\sigma}\big]$,
$$
\sum_{r \in {\mathbb Z} \setminus \{0, \pm 1\}} |u + r|^{-\mu} \le \frac{2}{\mu - 1}\, \big(1 -
\frac{1}{2\sigma}\big)^{1-\mu}\,.
$$
For all $n \in I_N$, we sustain
$$
|{\hat \varphi}(n \pm N_1)| = \frac{m}{N_1}\,\big|{\hat \varphi}_0 \big(\frac{m n}{N_1} \pm m \big)\big| \le \frac{m}{N_1}\,c_1\,,
$$
since
$$
\big|\frac{m n}{N_1} \pm m \big| \in \big[m \big(1 - \frac{1}{2\sigma}\big),\,m \big(1 + \frac{1}{2\sigma}\big)\big]\,.
$$
Thus we estimate for each $n \in I_N$,
\begin{eqnarray*}
\sum_{r \in {\mathbb Z} \setminus \{0\}} |{\hat \varphi}(n + r N_1)| &\le& \frac{m}{N_1}\,\Big[\big|{\hat \varphi}_0\big(\frac{m n}{N_1} - m \big)\big| + \big|{\hat \varphi}_0\big(\frac{m n}{N_1} + m \big)\big|\\
& & + \,\sum_{k \in {\mathbb Z} \setminus \{0, \pm 1\}} \big|{\hat \varphi}_0\big(\frac{m n}{N_1} + m r\big)\big|\Big]\\
&\le& \frac{m}{N_1}\, \Big[2 c_1 + \frac{c_2}{m^{\mu}}\,\sum_{r \in {\mathbb Z} \setminus \{0,\pm 1\}} \big|\frac{n}{N_1} + r \big|^{- \mu}\Big]\\
&\le& \frac{m}{N_1}\,\Big[2 c_1 + \frac{2 c_2}{(\mu - 1) \,m^{\mu}}\, \big(1 -\frac{1}{2\sigma}\big)^{1 -\mu}\Big]\,.
\end{eqnarray*}
Now we determine the minimum of all positive values
$$
{\hat \varphi}(n) = \frac{m}{N_1}\, {\hat \varphi}_0\big(\frac{m n}{N_1}\big)\,, \quad n \in I_N\,.
$$
Since $\frac{m\,|n|}{N_1} \le \frac{m}{2 \sigma}$ for all $n \in I_N$, we obtain
$$
\min_{n \in I_N} {\hat \varphi}(n) = \frac{m}{N_1}\,\min_{n \in I_N} {\hat \varphi}_0\big(\frac{m n}{N_1}\big)
= \frac{m}{N_1}\,{\hat \varphi}_0\big(\frac{m}{2 \sigma}\big) = {\hat \varphi}\big(\frac{N}{2}\big) > 0\,.
$$
Thus we see that the constant $e_{\sigma,N}(\varphi)$ can be estimated by an upper bound which depends
on $m$ and $\sigma$, but does not depend on $N$. We obtain
\begin{eqnarray*}
e_{\sigma,N}(\varphi) &\le& \frac{1}{{\hat \varphi}(N/2)}\,\max_{n \in I_N} \sum_{r \in {\mathbb Z} \setminus \{0\}} |{\hat \varphi}(n + r N_1)|\\
&\le& \frac{1}{{\hat \varphi}_0\big(\frac{m}{2\sigma}\big)}\,\Big[2 c_1 +
\frac{2 c_2}{(\mu - 1)\,m^{\mu}}\, \big(1 -\frac{1}{2\sigma}\big)^{1-\mu}\Big]\,.
\end{eqnarray*}
Consequently, the $C(\mathbb T)$-error constant $e_{\sigma}(\varphi)$ has the upper bound \eqref{eq:boundesigmascaledvarphi}. \qedsymbol

\subsection{B-spline window function}

We start with the popular B-spline window function (see \cite{Bey95,St97}).
Assume that $N\in 2\,\mathbb N$ and $\sigma > 1$ with $N_1 = \sigma N \in 2\,\mathbb N$ are given.
We consider the \emph{$\mathrm{B}$-spline window function}
\begin{equation}
\label{eq:B-spline}
\varphi_{\mathrm B} (x) := \frac{1}{M_{2m}(0)}\, M_{2m}(N_1 x)\,,
\end{equation}
where $M_{2m}$ denotes the centered cardinal B-spline of even order $2m$ with $m\in \mathbb N$. For $m=1$, we obtain
the \emph{triangular window function}.
Using the three-term recursion
\begin{equation}
\label{eq:Mm(x)}
M_m(x) = \frac{x + \frac{m}{2}}{m-1}\,M_{m-1}\big(x + \frac{1}{2}\big) + \frac{\frac{m}{2}-x}{m-1}\,M_{m-1}\big(x - \frac{1}{2}\big)\,, \quad m=2,\,3,\,\ldots\,,
\end{equation}
with
$$
M_1(x) := \left\{ \begin{array}{ll} 1 & \quad x \in (-\frac{1}{2},\,\frac{1}{2})\,,\\ [1ex]
\frac{1}{2} & \quad x\in \{-\frac{1}{2},\,\frac{1}{2}\}\,,\\ [1ex]
0 & \quad x\in \mathbb R \setminus [-\frac{1}{2},\,\frac{1}{2}]\,,
\end{array} \right.
$$
we find
$$
M_2(0) =1\,,\quad M_4(0) = \frac{2}{3}\,,\quad M_6(0) = \frac{11}{20}\,,\quad M_8(0) = \frac{151}{315}\,.
$$
{  Note that $M_{2m}(0) > 0$ for all $m \in \mathbb N$, since it holds $M_{2m}(x) > 0$ for each $x \in (-m,\,m)$.}
As known, the Fourier transform of \eqref{eq:B-spline} {  (see \cite[p. 452]{PlPoStTa18})} has the form
$$
{\hat \varphi}_{\mathrm B}(v) = \frac{1}{N_1\,M_{2m}(0)}\,\big(\mathrm{sinc}\frac{\pi v}{N_1}\big)^{2m}\,, \quad v \in \mathbb R\,,
$$
where
$$
\mathrm{sinc}\, x := \left\{\begin{array}{ll}
\sin x/x & \quad x \in \mathbb R \setminus \{0\}\,,\\
1 & \quad x=0\,.
\end{array}\right.
$$

If ${\tilde \varphi}_{\mathrm B}$ is the 1-periodization \eqref{eq:tildevarphi} of $\varphi_{\mathrm B}$, then the Fourier coefficients of $\tilde \varphi_{\mathrm B}$ read as follows
$$
c_k(\tilde \varphi_{\mathrm B}) = \int_0^1 {\tilde \varphi}_{\mathrm B}(t)\, {\mathrm e}^{- 2\pi {\mathrm i}\,k x}\, {\mathrm d}t = {\hat \varphi}_{\mathrm B}(k) = \frac{1}{N_1\,M_{2m}(0)}\,\big(\mathrm{sinc}\frac{\pi k}{N_1}\big)^{2m} \ge 0\,.
$$
Note that $c_k(\tilde \varphi_{\mathrm B}) > 0$ for all $k \in I_N$.

By \eqref{eq:s-fnew} we see that
\begin{equation*}
\|s - f\|_{C(\mathbb T)} \le  \sum_{n\in I_N} |c_n(f)| \,\sum_{r \in {\mathbb Z}\setminus \{0\}}  \frac{c_{n+r N_1}(\tilde \varphi_{\mathrm B})}{c_n(\tilde \varphi_{\mathrm B})}\,.
\end{equation*}

Now we estimate the $C(\mathbb T)$-error constants for NFFT. Applying the special structure of the Fourier coefficients $c_k({\tilde \varphi}_{\mathrm B})$ and Lemma \ref{Lemma3.7}, we obtain a good upper  bound \eqref{eq:boundesimple} of the $C(\mathbb T)$-error constant by this method.
For a proof of the following result see \cite{St97}.

\begin{Theorem}
Let $\sigma > 1$, $m \in \mathbb N \setminus \{1\}$, $N\in 2\,\mathbb N$, and $N_1 = \sigma N \in 2\,\mathbb N$. Further $m \in \mathbb N$ with $2m \ll N_1$ is given.
Then the $C(\mathbb T)$-error constant \eqref{eq:convenient1} of the $\mathrm{B}$-spline window function \eqref{eq:B-spline}   can be estimated by
$$
e_{\sigma}({\varphi}_{\mathrm B}) \le \frac{4m}{2m-1}\,(2\sigma -1)^{-2m}\,,
$$
i.e., the $\mathrm{B}$-spline window function \eqref{eq:B-spline} is convenient for $\mathrm{NFFT}$.
\end{Theorem}

\subsection{Modified B-spline window function}

Let $\sigma \ge 1$, $m\in \mathbb N \setminus \{1\}$, $N\in 2\,\mathbb N$, and $N_1 = \sigma N$ be given.
The approach to the B-spline window function \eqref{eq:B-spline} can be generalized to the \emph{modified $\mathrm{B}$-spline window function} (see \cite{Ne14})
\begin{equation}
\label{eq:modB-spline}
\varphi_{\mathrm{mB}} (x) := \frac{1}{M_{2b}(0)}\,M_{2b}(\frac{N_1 b}{m}\, x)\,,
\end{equation}
where $M_{2b}$ denotes the centered cardinal B-spline of order $2b\in \mathbb N \setminus \{1, \, 2\}$, i.e., $b \in \{\frac{3}{2},\,2,\,\frac{5}{2},\,\ldots\}$. Assume that $N_1 b \in 2\,\mathbb N$ and that $m \in \mathbb N$ fulfills the conditions $m < 2\sigma b$ and $b\not= m$.
Using the three-term recursion \eqref{eq:Mm(x)}, we find
$
M_3(0) =\frac{3}{4}$,  $M_5(0) = \frac{115}{192}\,.
$
Obviously, it holds $\varphi_{\mathrm{mB}} \in \Phi_{m, N_1}$, where the Fourier transform of \eqref{eq:modB-spline} {  (see \cite[p. 452]{PlPoStTa18})} reads as follows
$$
{\hat \varphi}_{\mathrm{mB}}(v) = \frac{m}{N_1 b}\,\big(\mathrm{sinc}\frac{m \pi v}{N_1 b }\big)^{2b}\,, \quad v \in \mathbb R \,.
$$
\medskip

If $\tilde \varphi_{\mathrm{mB}}$ denotes the 1-periodization of \eqref{eq:modB-spline}, then the Fourier coefficients of $\tilde \varphi_{\mathrm{mB}}$ have the following form
$$
c_k(\tilde \varphi_{\mathrm{mB}}) = \int_0^1 {\tilde \varphi}_{\mathrm{mB}}(t)\, {\mathrm e}^{- 2\pi {\mathrm i}\,k x}\, {\mathrm d}t = {\hat \varphi}_{\mathrm{mB}}(k) = \frac{m}{N_1 b}\,\big(\mathrm{sinc}\frac{m \pi k}{N_1 b }\big)^{2b} \,.
$$
Note that $c_k({\tilde \varphi}_{\mathrm{mB}}) > 0$ for all $k \in I_N$.
Let $f$
be an arbitrary 1-periodic trigonometric polynomial \eqref{eq:f(x)} which we approximate by the 1-periodic function
$$
s(x) := \sum_{\ell \in I_{N_1 b}} g_{\ell}\,{\tilde \varphi_{\mathrm{mB}}}(x - \frac{\ell}{N_1 b})
$$
with conveniently chosen coefficients $g_{\ell} \in \mathbb C$. Then $s$ possesses the Fourier expansion
$$
s(x) = \sum_{k\in \mathbb Z} c_k(s)\, {\mathrm e}^{2\pi {\mathrm i}\,k x}
$$
with the Fourier coefficients
$$
c_k(s) = \int_0^1 s(t)\, {\mathrm e}^{- 2\pi {\mathrm i}\,k x}\, {\mathrm d}t = {\hat g}_k\, c_k(\tilde \varphi_{\mathrm{mB}})\,,
$$
where
$$
{\hat g}_k = \sum_{\ell \in I_{N_1 b}} g_{\ell}\, {\mathrm e}^{- 2\pi {\mathrm i}\,k \ell/(N_1 b)}\,.
$$
Thus the vector $({\hat g}_k)_{k^\in I_{N_1 b}}$ is equal to the DFT of length $N_1 b$ of the vector $(g_{\ell})_{\ell\in I_{N_1 b}}$ and we have
${\hat g}_{k + r N_1 b} = {\hat g}_k$ for all $k\in I_{N_1 b}$ and $r\in \mathbb Z$.
In order to approximate $f$ by $s$, we choose
$$
{\hat g}_k = \left\{ \begin{array}{ll} \frac{c_k(f)}{c_k(\tilde \varphi_{\mathrm{mB}})} & \quad k \in I_N\,,\\ [1ex]
0 & \quad k \in I_{N_1 b} \setminus I_N\,.
\end{array} \right.
$$
Thus we see that $c_k(s) = c_k(f)$ for all $k \in I_N$ and $c_k(s) = 0$ for all $k \in I_{N_1 b} \setminus I_N$. Substituting $k = n + r N_1 b$ with $n\in I_{N_1 b}$
and $r \in {\mathbb Z}$, we obtain
\begin{eqnarray*}
s(x) - f(x) &=& \sum_{k \in \mathbb Z \setminus I_{N_1 b}} c_k(s)\,{\mathrm e}^{2\pi {\mathrm i}\,k x} = \sum_{n\in I_{N_1 b}} \sum_{r \in {\mathbb Z}\setminus \{0\}} c_{n+r N_1 b}(s)\,{\mathrm e}^{2\pi {\mathrm i}\,(n + r N_1 b)\,x}\\
&=& \sum_{n\in I_N} \sum_{r \in {\mathbb Z}\setminus \{0\}} c_n(f)\, \frac{1}{c_n(\tilde \varphi_{\mathrm{mB}})}\,c_{n+r N_1 b}(\tilde \varphi_{\mathrm{mB}})\,{\mathrm e}^{2\pi {\mathrm i}\,(n + r N_1 b)\,x}
\end{eqnarray*}
such that
\begin{equation*}
|s(x) - f(x)| \le  \sum_{n\in I_N} |c_n(f)| \,\sum_{r \in {\mathbb Z}\setminus \{0\}}  \frac{|c_{n+r N_1 b}(\tilde \varphi_{\mathrm{mB}})|}{c_n(\tilde \varphi_{\mathrm{mB}})}\,.
\end{equation*}
Then it follows by H\"{o}lder's inequality that
\begin{eqnarray*}
\|s -f\|_{C(\mathbb T)}
&\le& \Big( \sum_{n\in I_N} |c_n(f)|\Big)\, \max_{n\in I_N}\sum_{r \in {\mathbb Z}\setminus \{0\}}  \frac{|c_{n+r N_1 b}(\tilde \varphi_{\mathrm{mB}})|}{c_n(\tilde \varphi_{\mathrm{mB}})}\\[1ex]
&=& \|f\|_{A(\mathbb T)} \, \max_{n\in I_N}\sum_{r \in {\mathbb Z}\setminus \{0\}}  \frac{|c_{n+r N_1 b}(\tilde \varphi_{\mathrm{mB}})|}{c_n(\tilde \varphi_{\mathrm{mB}})}\,.
\end{eqnarray*}

\begin{Theorem}
Let $\sigma \ge 1$, $N\in 2\,\mathbb N$, $2 b \in \mathbb N \setminus \{1,\,2\}$, and $N_1 b\in 2\,\mathbb N$ with $N_1 = \sigma N$. Further let $m \in \mathbb N$ with $2 b > m$ and $b\not= m$ be given.

Then the $C(\mathbb T)$-error constant of $\mathrm{NFFT}$ with the modified $\mathrm{B}$-spline window function \eqref{eq:modB-spline} can be estimated by
$$
e_{\sigma b}(\varphi_{\mathrm{mB}}) \le \frac{4 b}{2b - 1}\,(2 \sigma b - 1)^{-2b} \,,
$$
i.e., the modified  $\mathrm{B}$-spline window function \eqref{eq:modB-spline} is convenient for  $\mathrm{NFFT}$.
\end{Theorem}

\emph{Proof.} Now we estimate
$$
\max_{n\in I_N}\,\sum_{r \in {\mathbb Z}\setminus \{0\}}  \frac{|c_{n+r N_1 b}(\tilde \varphi_{\mathrm{mB}})|}{c_n(\tilde \varphi_{\mathrm{mB}})}\,.
$$
For $n=0$ and $r \in \mathbb Z \setminus \{0\}$ we have
$$
\frac{|c_{r N_1 b}(\tilde \varphi_{\mathrm{mB}})|}{1 \cdot c_0(\tilde \varphi_{\mathrm{mB}})} = 0\,.
$$
For $n \in I_N \setminus \{0\}$ and  $r \in \mathbb Z \setminus \{0\}$ we obtain
\begin{eqnarray*}
\frac{|c_{n+r N_1 b}(\tilde \varphi_{\mathrm{mB}})|}{c_n(\tilde \varphi_{\mathrm{mB}})} &=& \frac{|n|^{2b}}{|n + r N_1 b|^{2b}} = \big(\frac{|n|}{N_1 b}\big)^{2b}\,\big| \frac{n}{N_1 b} + r \big|^{-2b}\\
&\le& \big(\frac{1}{2\sigma b}\big)^{2b} \,|\frac{n}{N_1 b} +  r|^{-2b}\,,
\end{eqnarray*}
since for $n \in I_N$ it holds
$$
\frac{|n|}{N_1 b} \le \frac{1}{2 \sigma b} \le \frac{1}{3}\,.
$$
Using Lemma \ref{Lemma3.7}, we conclude that
\begin{eqnarray*}
\sum_{r \in {\mathbb Z}\setminus \{0\}}  \frac{|c_{n+r N_1 b}(\tilde \varphi_{\mathrm{mB}})|}{c_n(\tilde \varphi_{\mathrm{mB}})} &\le& \frac{1}{(2 \sigma b)^{2b}} \sum_{r\in {\mathbb Z}\setminus \{0\}}
\big| \frac{n}{N_1 b} + r \big|^{-2b}\\
&\le& \frac{1}{(2 \sigma b)^{2b}}\Big[2\,\big(1 - \frac{|n|}{N_1 b}\big)^{-2b} + \sum_{r\in {\mathbb Z}\setminus \{0, \pm 1\}}\big| \frac{n}{N_1 b} + r \big|^{-2b}\Big]\\
&\le& \frac{1}{(2 \sigma b)^{2b}}\Big[2\,\big(1 - \frac{1}{2\sigma b}\big)^{-2b} + \frac{2}{2b-1}\,\big(1 - \frac{1}{2\sigma b}\big)^{1-2b}\Big]\\
&\le&  \frac{4b}{(2b - 1)\,(2 \sigma b)^{2b}}\,\big( 1 - \frac{1}{2 \sigma b})^{-2b} = \frac{4 b}{2b - 1}\,(2 \sigma b - 1)^{-2b}
\end{eqnarray*}
with $2 \sigma b \ge 3$. \qedsymbol

\subsection{Algebraic window function}

For fixed shape parameter  $\beta = 3 m$ with $m \in \mathbb N \setminus \{1\}$ and for an oversampling factor $\sigma > \frac{\pi}{3}$, we consider the \emph{algebraic window function}
\begin{equation}
\label{eq:algwindow}
\varphi_{\mathrm{alg}}(x) := \left\{ \begin{array}{ll}
\big(1-\frac{(N_1 x)^2}{m^2}\big)^{\beta-1/2} & \quad x \in \big[- \frac{m}{N_1},\,\frac{m}{N_1}\big]\,, \\ [1ex]
0 & \quad x \in \mathbb R \setminus \big[-\frac{m}{N_1},\,\frac{m}{N_1}\big]\,.
\end{array}\right.
\end{equation}

\begin{Theorem}
\label{Thm:esigmaNpol}
Let $N \in 2 \mathbb N$ and $\sigma > \frac{\pi}{3}$, where $N_1 = \sigma N \in 2 \mathbb N$. Further let $m \in \mathbb N \setminus \{1\}$  and $\beta = 3m$ be given.

Then the $C(\mathbb T)$-error constant of the algebraic window function \eqref{eq:algwindow} can be estimated by
\begin{equation}
\label{eq:esigmavarphialgest}
e_{\sigma}({\varphi}_{\mathrm{alg}}) \le \frac{3 \sqrt{\sigma}}{\sqrt{\pi m}\,J_{3m}\big(\frac{\pi m}{\sigma}\big)}\,\Big[1 + \frac{2 \sigma - 1}{(6m-1)\sigma }\Big]\,(2 \sigma - 1)^{-3m-1/2}\,,
\end{equation}
i.e., the algebraic window function \eqref{eq:algwindow} is convenient for ${\mathrm{NFFT}}$.
\end{Theorem}

\emph{Proof.} We apply Theorem \ref{Lemma:esigmaNvarphibounded} and introduce the unscaled algebraic window function
$$
\varphi_{0,\mathrm{alg}}(x) := \left\{ \begin{array}{ll}
(1-x^2)^{\beta-1/2} & \quad x \in [- 1,\,1]\,, \\
0 & \quad x \in \mathbb R \setminus \big[-1,\,1]\,.
\end{array}\right.
$$
Using  \cite[p.~8]{Oberh90}, we determine the corresponding Fourier transform
\begin{eqnarray}
\label{eq:FTvarphi0alg}
{\hat \varphi}_{0, \mathrm{alg}}(v)
&=& \int_{-1}^1 (1-x^2)^{\beta-1/2}\,{\mathrm e}^{-2\pi {\mathrm i}\, v x}\, {\mathrm d}x = 2\,\int_0^1 (1-x^2)^{\beta-1/2}\,\cos(2\pi\, v x)\, {\mathrm d}x \nonumber \\ [1ex]
&=& \frac{\pi \,(2\beta)!}{4^{\beta}\,\beta!} \left\{\begin{array}{ll}
(\pi v)^{-\beta}\,J_{\beta}(2\pi  v)& \quad v \in {\mathbb R}\setminus \{0\}\,,\\
\frac{1}{\beta!} & \quad v=0\,.
\end{array} \right.
\end{eqnarray}
Thus ${\hat \varphi}_{0, \mathrm{alg}}(v)$ for $v>0$ is a multiple of the function
$(\pi v)^{-\beta}\,J_{\beta}(2 \pi v)$
which can be represented as infinite product
\begin{equation}
\label{eq:infproduct}
\frac{1}{\beta!}\,\prod_{s=1}^{\infty} \big(1 - \frac{(2\pi v)^2}{j_{\beta,s}^2}\big)\,,
\end{equation}
where $j_{\beta,s}$ denotes the $s$th positive zero of $J_{\beta}$ (see \cite[p.~370]{abst}).
Note that
\begin{eqnarray*}
j_{6,1} &=& 9.936109\ldots\,, \quad j_{9,1} = 13.354300\ldots\,, \quad j_{12,1} = 16.698249\ldots\,,\\
j_{15,1} &=& 19.994430\ldots\,, \quad j_{18,1} = 23.256776\ldots\,, \quad j_{21,1} = 26.493647\ldots\,.
\end{eqnarray*}
For $\beta = 3m$ it holds $j_{\beta,1} > 3m + \pi - \frac{1}{2}$ (see \cite{IfSi85}). Hence by $\sigma > \frac{\pi}{3}$ we get
$$
\frac{2 \pi m}{2 \sigma\,j_{\beta,1}} < \frac{ \frac{\pi}{\sigma}\,m}{3m + \pi -\frac{1}{2}} < \frac{3m}{3m + \pi -\frac{1}{2}} <1\,.
$$
{  Thus each factor of the infinite product \eqref{eq:infproduct} is positive and decreasing for  $v \in \big[0,\, \frac{m}{2\sigma}\big]$. Hence}
by \eqref{eq:FTvarphi0alg} and \eqref{eq:infproduct}, the Fourier transform ${\hat \varphi}_{0, \mathrm{alg}}(v)$ is positive and decreasing for $v \in \big[0,\, \frac{m}{2\sigma}\big]$ and it holds
\begin{equation}
\label{eq:FTvarphi0alg(m)}
{\hat \varphi}_{0,\mathrm{alg}}\big(\frac{m}{2\sigma}\big) = \frac{(2\beta)!\, \pi}{4^{\beta}\,\beta!}\,\big(\frac{2 \sigma}{\pi m}\big)^{\beta}\, J_{\beta}\big(\frac{\pi m}{\sigma}\big)\,.
\end{equation}
By \eqref{eq:|J3m(x)|} we know that for $|v|\ge \frac{m}{2}$ it holds
\begin{equation}
\label{eq:|Jbeta(2piv)|}
|J_{\beta}(2 \pi v)| \le \frac{3}{2 \,\sqrt{2 \pi |v|}}\,.
\end{equation}
For $|v| \ge m\,\big(1 - \frac{1}{2 \sigma}\big)$ and $\sigma > \frac{\pi}{3}$ it holds
$$
|v| \ge m\,\big(1 - \frac{1}{2 \sigma}\big) > \frac{m}{2}\,.
$$
Using \eqref{eq:FTvarphi0alg} and \eqref{eq:|Jbeta(2piv)|}, we obtain for $|v| \ge m\,\big(1 - \frac{1}{2 \sigma}\big)$ the estimate
\begin{equation}
\label{eq:FTvarphi0algestim}
|{\hat \varphi}_{0,\mathrm{alg}}(v)| \le \frac{3\,(2  \beta)!}{2^{3/2}\,4^{\beta}\,\beta!\,\pi^{\beta-1/2} }\, |v|^{-\beta-1/2}\,.
\end{equation}
Applying Theorem \ref{Lemma:esigmaNvarphibounded}, we obtain
$$
e_{\sigma}(\varphi) \le \frac{1}{{\hat \varphi}_0\big(\frac{m}{2 \sigma}\big)} \Big[2 c_1 + \frac{2 c_2}{(\mu -1)\,m^{\mu}}\big(1 - \frac{1}{2\sigma}\big)^{1-\mu}\Big]\,,
$$
where by \eqref{eq:FTvarphi0algestim} we have $\mu= \beta + \frac{1}{2} = 3m + \frac{1}{2}$ and
\begin{eqnarray*}
c_1 &=& \frac{3\,(2  \beta)!}{2^{3/2}\,4^{\beta}\,\beta!\,\pi^{\beta-1/2} }\, m^{- \beta -1/2}\, \big(1 - \frac{1}{2\sigma}\big)^{1- \mu}\,,\\
c_2 &=& \frac{3\,(2  \beta)!}{2^{3/2}\,4^{\beta}\,\beta!\,\pi^{\beta-1/2} }\,.
\end{eqnarray*}
Using \eqref{eq:FTvarphi0alg(m)}, we get the inequality \eqref{eq:esigmavarphialgest}. \qedsymbol

\subsection{Bessel window function}

Let $N \in 2\, \mathbb N$, $m \in \mathbb N \setminus \{1\}$, and $\sigma \in \big[\frac{5}{4},\,2\big]$ be given. For fixed shape parameter
\begin{equation}
\label{eq:shapeparam}
\beta := 2 \pi m\, \big(1 - \frac{1}{2\sigma}\big)\,,
\end{equation}
we consider the new \emph{Bessel window function}
\begin{equation}
\label{eq:Besselwindow}
\varphi_{\mathrm{Bessel}}(x) := \left\{\begin{array}{ll} \frac{1}{I_2(\beta)}\,\big(1 - \frac{(N_1x)^2}{m^2}\big)\,I_2\big(\beta\,\sqrt{1 - \frac{(N_1x)^2}{m^2}}\big) & \quad x\in \big[-\frac{m}{N_1},\,\frac{m}{N_1}\big]\,,\\ [1ex]
0 & \quad x \in {\mathbb R} \setminus \big[-\frac{m}{N_1},\,\frac{m}{N_1}\big]
\end{array} \right.
\end{equation}
with $N_1 = \sigma N \in 2 \mathbb N$. Obviously, the Bessel window function is continuously differentiable and compactly supported.

\begin{Theorem}
\label{Thm:esigmaNBessel}
Let $N \in 2 \mathbb N$,  $m \in \mathbb N \setminus \{1\}$, and $\sigma \in \big[\frac{5}{4},\,2\big]$ be given, where $N_1 = \sigma N \in 2 \mathbb N$ and $2m \ll N_1$.

Then the $C(\mathbb T)$-error constant of the Bessel window function \eqref{eq:Besselwindow} with the shape parameter \eqref{eq:shapeparam} can be estimated by
$$
e_{\sigma}({\varphi}_{\mathrm{Bessel}}) \le \big(50\,m^3 + 7\big) \,{\mathrm e}^{- 2 \pi m\, \sqrt{1 - 1/\sigma}}\,.
$$
i.e., the Bessel window function \eqref{eq:Besselwindow} is convenient for ${\mathrm{NFFT}}$.
\end{Theorem}

\emph{Proof.} We apply Theorem \ref{Lemma:esigmaNvarphibounded} and introduce the unscaled Bessel window function
$$
\varphi_{0,\mathrm{Bessel}}(x) := \left\{\begin{array}{ll} \frac{1}{I_2(\beta)}\,(1 - x^2)\,I_2\big(\beta\,\sqrt{1 - x^2}\big) & \quad x\in [-1,\,1]\,,\\
0 & \quad x \in {\mathbb R} \setminus [-1,\,1]\,.
\end{array} \right.
$$
We determine the even Fourier transform
\begin{eqnarray*}
{\hat \varphi}_{0,\mathrm{Bessel}}(v) &=& \int_{\mathbb R} \varphi_{0,\mathrm{Bessel}}(x)\,{\mathrm e}^{-2\pi {\mathrm i}\,v x}\,{\mathrm d}x \\ [1ex]
&=& \frac{2}{I_2(\beta)}\,\int_0^1 (1 - x^2)\,I_2\big(\beta\,\sqrt{1 - x^2}\big)\,\cos (2\pi\, v x)\,{\mathrm d}x\,.
\end{eqnarray*}
By \cite[p.~96]{Oberh90}, this Fourier transform reads as follows
\begin{equation}
\label{eq:FTvarphi0Bessel}
{\hat \varphi}_{0,\mathrm{Bessel}}(v) = \frac{2\,\beta^2}{I_2(\beta)}\,\left\{ \begin{array}{ll}
\sqrt{\frac{\pi}{2}}\,(\beta^2 - 4 \pi^2 v^2)^{-5/4}\,I_{5/2}\big(\sqrt{\beta^2 - 4 \pi^2 v^2}\big) & \quad |v| < m\,\big(1 - \frac{1}{2 \sigma}\big)\,,\\ [1ex]
\frac{1}{15} & \quad v = \pm m\,\big(1 - \frac{1}{2 \sigma}\big)\,,\\ [1ex]
\sqrt{\frac{\pi}{2}}\,(4 \pi^2 v^2 -\beta^2)^{-5/4}\,J_{5/2}\big(\sqrt{4 \pi^2 v^2 - \beta^2}\big) & \quad |v| > m\,\big(1 - \frac{1}{2 \sigma}\big)\,.
\end{array} \right.
\end{equation}
Introducing the \emph{spherical Bessel function} (see \cite[pp.~437--438]{abst})
$$
j_2(x) := \sqrt{\frac{\pi}{2x}}\, J_{5/2}(x) = \big(\frac{3}{x^3} - \frac{1}{x}\big)\,\sin x - \frac{3}{x^2}\,\cos x\,, \quad x > 0; \; j_2(0) := 0\,,
$$
and the \emph{modified spherical Bessel function} (see \cite[p.~443]{abst})
$$
i_2(x) := \sqrt{\frac{\pi}{2x}}\, I_{5/2}(x) = \big(\frac{3}{x^3} + \frac{1}{x}\big)\,\sinh x - \frac{3}{x^2}\,\cosh x\,, \quad x > 0; \; i_2(0) := 0\,,
$$
we obtain
\begin{equation}
\label{eq:hatvarphi0Bessel}
{\hat \varphi}_{0,\mathrm{Bessel}}(v) = \frac{2 \beta^2}{I_2(\beta)}\,\left\{ \begin{array}{ll}
(\beta^2 - 4 \pi^2 v^2)^{-1}\,i_2\big(\sqrt{\beta^2 - 4 \pi^2 v^2}\big) & \quad |v| <  m\,\big(1 - \frac{1}{2 \sigma}\big)\,,\\ [1ex]
\frac{1}{15} & \quad v = \pm  m\,\big(1 - \frac{1}{2 \sigma}\big)\,,\\ [1ex]
(4 \pi^2 v^2 - \beta^2)^{-1}\,j_2\big(\sqrt{4 \pi^2 v^2 - \beta^2}\big) & \quad |v| >  m\,\big(1 - \frac{1}{2 \sigma}\big)\,.
\end{array} \right.
\end{equation}
Using the power series expansion of the modified spherical Bessel function $i_2$ (see \cite[p.~443]{abst}), for $|v| <  m\,\big(1 - \frac{1}{2 \sigma}\big)$ we receive
$$
(\beta^2 - 4 \pi^2 v^2)^{-1}\,i_2\big(\sqrt{\beta^2 - 4 \pi^2 v^2}\big) = \sum_{k=0}^{\infty} \frac{1}{2^k\,k!\,(2k+5)!!}\,(\beta^2 - 4 \pi^2 v^2)^k\,.
$$
Hence ${\hat \varphi}_{0,\mathrm{Bessel}}(v)$ is positive and decreasing for $v \in \big[0,\, m \big(1 - \frac{1}{2\sigma}\big) \big)$. Since $\frac{m}{2 \sigma} < m \big(1 - \frac{1}{2 \sigma}\big)$ for $ \sigma \ge \frac{5}{4}$, we conclude
$$
{\hat \varphi}_{0,\mathrm{Bessel}}\big(\frac{m}{2\sigma}\big) = \frac{2}{I_2(\beta)}\,\big(1 - \frac{1}{2\sigma}\big)^{2}\,\big(1 - \frac{1}{\sigma}\big)^{-1}\,i_2\big(2 \pi m\,\sqrt{1- \frac{1}{\sigma}}\big)\,.
$$
For $m\ge 2$ and $\sigma \ge \frac{5}{4}$, we maintain
$$
2 \pi m\,\sqrt{1- \frac{1}{\sigma}} \ge 4 \pi\,\sqrt{1 - \frac{1}{\sigma}} \ge x_0 := \frac{4 \pi}{\sqrt 5}\,.
$$
Applying the inequality \eqref{eq:Imu(x)estim} with $\mu = \frac{5}{2}$, for $x \ge x_0$ we obtain
{ 
$$
i_2(x) \ge  x_0 \,{\mathrm e}^{-x_0}\,i_2(x_0)\, x^{-1}\,{\mathrm e}^{x} > 0.280573 \, x^{-1}\,{\mathrm e}^{x} > \frac{1}{4}\, x^{-1}\,{\mathrm e}^{x}\,.
$$
}
Hence for $x = 2 \pi m \sqrt{1 - \frac{1}{\sigma}}$ it follows that
$$
i_2\big(2 \pi m \sqrt{1 - \frac{1}{\sigma}}\big) > \frac{1}{8 \pi m}\,\big(1 - \frac{1}{\sigma}\big)^{-1/2}\,{\mathrm e}^{2 \pi m \sqrt{1 - 1/\sigma}}\,.
$$
Thus we see that
\begin{equation}
\label{eq:hatvarphi0Bessel(m)}
{\hat \varphi}_{0,\mathrm{Bessel}}\big(\frac{m}{2\sigma}\big) \ge \frac{1}{4 \pi m\,I_2(\beta)}\, \big(1-\frac{1}{2 \sigma}\big)^2\,\big(1-\frac{1}{\sigma}\big)^{-3/2}\,{\mathrm e}^{2 \pi m\, \sqrt{1 - 1/\sigma}}\,.
\end{equation}
Now we estimate the Fourier transform ${\hat \varphi}_{0,\mathrm{Bessel}}(v)$ for $|v| \ge m\,\big(1 + \frac{1}{2 \sigma}\big)$.
Using the assumptions $m \ge 2$ and $\sigma \in \big[\frac{5}{4},\,2\big]$, for $|v| \ge m\,\big(1 + \frac{1}{2 \sigma}\big)$ we get
\begin{eqnarray*}
\sqrt{4 \pi^2 v^2 - \beta^2} &=& \sqrt{4 \pi^2 v^2 - 4 \pi^2 m^2\, \big(1 - \frac{1}{2 \sigma}\big)^2} \ge 2 \pi m\, \sqrt{\big(1 + \frac{1}{2 \sigma}\big)^2 - \big(1 - \frac{1}{2 \sigma}\big)^2}\\
&=& 2 \pi m \,\sqrt{\frac{2}{\sigma}} \ge 4\pi\,\sqrt{\frac{2}{\sigma}} \ge 4 \pi > 6\,.
\end{eqnarray*}
By \eqref{eq:|J5/2(x)|}, for all $x \ge 6$ it holds
$$
|J_{5/2}(x)| < \frac{1}{\sqrt x}\,.
$$
Thus for $|v| \ge m\, \big(1 + \frac{1}{2 \sigma}\big)$, we obtain by \eqref{eq:FTvarphi0Bessel} that
$$
|{\hat \varphi}_{0,\mathrm{Bessel}}(v)| \le \frac{\sqrt{2 \pi} \,\beta^2}{8 \pi^3\,I_2(\beta)}\,\big(v^2 - m^2 \big(1 - \frac{1}{2 \sigma}\big)^2\big)^{-3/2}\,.
$$
Since the function $g:\,\big[m\,\big(1 + \frac{1}{2\sigma}\big),\, \infty \big) \to \mathbb R$ defined by
$$
g(v):= v^3\, \big(v^2 - m^2 \big(1 - \frac{1}{2\sigma}\big)^2\big)^{-3/2}=\big(1 - \frac{m^2}{v^2}\,\big(1 - \frac{1}{2\sigma}\big)^2\big)^{-3/2}\,,
$$
is decreasing and {  bounded from above} by
$$
g\big(m \big(1 + \frac{1}{2 \sigma}\big) \big) = \frac{(2 \sigma + 1)^3}{ (8 \sigma)^{3/2}}\,,
$$
we receive for $|v| \ge m\, \big(1 + \frac{1}{2 \sigma}\big)$ the estimate
\begin{equation}
\label{eq:|hatvarphi0Bessel(v)|}
|{\hat \varphi}_{0,\mathrm{Bessel}}(v)| \le \frac{m^2\,(2 \sigma + 1)^3}{32 \sigma \sqrt{\pi \sigma} \, I_2(\beta)} \,\big(1 - \frac{1}{2 \sigma}\big)^2\, |v|^{-3}\,.
\end{equation}
Finally we show that for $|v| \ge m \,\big(1 - \frac{1}{2 \sigma}\big)$ it holds
\begin{equation}
\label{eq:|hatvarphi0Bessel(smallv)|}
|{\hat \varphi}_{0,\mathrm{Bessel}}(v)| \le \frac{2 \beta^2}{15\,I_2(\beta)}\,.
\end{equation}
By \eqref{eq:hatvarphi0Bessel} we have
$$
{\hat \varphi}_{0,\mathrm{Bessel}}\big(\pm m \big(1 - \frac{1}{2 \sigma}\big)\big) = \frac{2 \beta^2}{15\,I_2(\beta)}\,.
$$
For $|v| > m \big(1 - \frac{1}{2 \sigma}\big)$, it holds by \eqref{eq:hatvarphi0Bessel} that
$$
{\hat \varphi}_{0,\mathrm{Bessel}}(v) = \frac{2 \beta^2}{I_2(\beta)}\,(4 \pi^2 v^2 - \beta^2)^{-1}\,j_2\big(\sqrt{4 \pi^2 v^2 - \beta^2}\big)\,.
$$
{  By the definition of the spherical Bessel function $j_2$ it holds for $x> 0$ that
$$
|j_2(x)| = \sqrt{\frac{\pi}{2x}}\, |J_{5/2}(x)|\,.
$$
Then from \cite[p. 49]{Watson} it follows that for $x > 0$,
$$
|J_{5/2}(x)| \le \frac{1}{\Gamma(7/2)}\,\big(\frac{x}{2}\big)^{5/2}
$$
with $\Gamma(7/2) = \frac{15}{8}\,\sqrt \pi$ such that by $j_2(0) = 0$ we obtain the inequality
$$
|j_2(x)| \le \frac{1}{15}\,x^2\,, \quad x \ge 0\,,
$$
and hence \eqref{eq:|hatvarphi0Bessel(smallv)|}.}

Applying \eqref{eq:boundesigmascaledvarphi} together with  \eqref{eq:hatvarphi0Bessel(m)}, \eqref{eq:|hatvarphi0Bessel(v)|}, and \eqref{eq:|hatvarphi0Bessel(smallv)|}, we conclude that
$$
e_{\sigma}({\varphi}_{\mathrm{Bessel}}) \le \Big[\frac{(4 \pi)^3}{15}\,\big(1 - \frac{1}{\sigma}\big)^{3/2}\,m^3 + c(\sigma)\Big] \,{\mathrm e}^{- 2 \pi m\, \sqrt{1 - 1/\sigma}}\,.
$$
with
$$
c(\sigma) := \frac{\sqrt \pi\,(2 \sigma + 1)^3}{8 \, \sigma^{3/2}}\,\big(1 - \frac{1}{2\sigma}\big)^{-2}\, \big(1 - \frac{1}{\sigma}\big)^{3/2}\,.
$$
For $\sigma \in \big[\frac{5}{4}, \,2\big]$, it holds
$$
\frac{(4 \pi)^3}{15}\,\big(1 - \frac{1}{\sigma}\big)^{3/2} < 50\,, \quad c(\sigma) < 7\,.
$$
This completes the proof. \qedsymbol
\medskip

\subsection{sinh-type and related window functions}\label{Sec:window}

Let $N \in 2 \mathbb N$, $m \in \mathbb N \setminus \{1\}$, and $\sigma \in \big[\frac{5}{4},\,2\big]$ be given. For fixed shape parameter {  $\beta$ as stated in} \eqref{eq:shapeparam},
we consider the new {$\sinh$-\emph{type window function}
\begin{equation}
\label{eq:sinhwindow}
\varphi_{\sinh} (x) := \left\{ \begin{array}{ll}
\frac{1}{\sinh \beta}\,\sinh \big(\beta \sqrt{1 - \frac{(N_1 x)^2}{m^2}}\big)   &\quad  x \in \big[ - \frac{m}{N_1},\, \frac{m}{N_1}\big]\,,\\
0                   & \quad x \in \mathbb R \setminus \big[ - \frac{m}{N_1},\, \frac{m}{N_1}\big]
\end{array} \right.
\end{equation}
with $N_1 = \sigma N \in 2 \mathbb N$. This $\sinh$-type window function belongs to $\Phi_{m,N_1}$. Note that \eqref{eq:sinhwindow} is not piecewise continuously differentiable, since $\varphi_{\sinh}'(-\frac{m}{N_1}+0)=\infty$
and $\varphi_{\sinh}'(\frac{m}{N_1}-0)= -\infty$.
Note that the $\sinh$-type window function \eqref{eq:sinhwindow} can be computed much faster than the Bessel window function \eqref{eq:Besselwindow}.

\begin{Theorem}
\label{Thm:esigmaNsinh}
Let $N \in 2 \mathbb N$, $N \ge 8$, $m \in \mathbb N \setminus \{1\}$, and $\sigma \in \big[\frac{5}{4},\,2\big]$ be given, where $N_1 = \sigma N \in 2 \mathbb N$ and $2m \ll N_1$.
{  The shape parameter $\beta$ is given by \eqref{eq:shapeparam}.}

Then the $C(\mathbb T)$-error constant of the $\sinh$-type window function \eqref{eq:sinhwindow} with the shape parameter \eqref{eq:shapeparam} can be estimated by
$$
e_{\sigma}({\varphi}_{\sinh}) \le \big(24\,m^{3/2} + 3\big)\,{\mathrm e}^{- 2 \pi m\, \sqrt{1 - 1/\sigma}}\,,
$$
i.e., the $\sinh$-type window function \eqref{eq:sinhwindow} is convenient for ${\mathrm{NFFT}}$.
\end{Theorem}

{  A proof of Theorem \ref{Thm:esigmaNsinh} can be found in \cite{PoTa20}. The proof based mainly on the knowledge of the analytical Fourier transform. In  \cite{PoTa20} we consider in addition two related window functions,
namely the \emph{continuous} $\exp$-\emph{type window function}
\begin{equation*}
\varphi_{\mathrm{exp}}(x) := \left\{ \begin{array}{ll}
\frac{1}{{\mathrm e}^{\beta} - 1} \Big({\mathrm e}^{\beta \,\sqrt{1 - (N_1x/m)^2}} - 1\Big) & \quad x \in I\,, \\ [1ex]
0 & \quad x \in \mathbb R \setminus  I\,.
\end{array}\right.
\end{equation*}
as well as the \emph{continuous} $\cosh$-\emph{type window function}
\begin{equation*} 
\varphi_{\cosh}(x) := \left\{ \begin{array}{ll}
\frac{1}{\cosh \beta - 1}\,\Big(\cosh\left(\beta \sqrt{1- \left(\frac{N_1 x}{m}\right)^2 }\right) - 1\Big) & \quad x \in I\,, \\ [1ex]
0 & \quad x \in \mathbb R \setminus  I\,.
\end{array}\right.
\end{equation*}
The main drawback for the numerical analysis of the $\exp$-type/$\cosh$-type window function is the fact that an explicit
Fourier transform of this window function is unknown. Therefore we split the continuous $\exp$-type/$\cosh$-type window function
into a sum $\psi + \rho$, where the Fourier transform of the compactly supported function $\psi$ is explicitly known and where the compactly supported function $\rho$ has small magnitude.  Note that $\varphi_{\mathrm{exp}}$ and was first suggested in \cite{BaMaKl18,Bar} and a discontinuous version of $\varphi_{\cosh}$ was suggested in \cite[Remark 13]{BaMaKl18}.
}

\medskip

\subsection{Modified cosh-type and related window functions}\label{Sec:modwindow}

{ 
For $\sigma \in \big[\frac{5}{4},\,2 \big]$  and $m \in \mathbb N \setminus \{1\}$, we consider the \emph{modified} $\cosh$-\emph{type window function}
\begin{equation}
\label{eq:mcoshwindow}
\varphi_{\mathrm{mcosh}}(x) := \varphi_{0,\mathrm{mcosh}}\big(\frac{N_1 x}{m}\big)\,, \quad x \in \mathbb R\,,
\end{equation}
where it holds
\begin{equation*} 
\varphi_{0,\mathrm{mcosh}} (x) := \left\{ \begin{array}{ll}
\frac{1}{\cosh \beta - 1}\,\frac{\cosh \big(\beta \sqrt{1 - x^2}\big)- 1}{\sqrt{1 - x^2}}    &\quad  x \in ( - 1,\, 1)\,,\\
0                   & \quad x \in \mathbb R \setminus ( - 1,\, 1)
\end{array} \right.
\end{equation*}
with the shape parameter $\beta$ as stated in \eqref{eq:shapeparam}. By \cite[p.~6 and p.~38]{Oberh90}, its Fourier transform reads as follows
\begin{eqnarray}
\label{eq:hatomegamcosh(v)}
& &{\hat \varphi}_{0,\mathrm{mcosh}}(v) = \frac{2}{\cosh \beta - 1}\,\int_0^1 \frac{\cosh\big(\beta \sqrt{1 - x^2}\big)- 1\big)}{\sqrt{1 - x^2}}\,\cos(2\pi v x)\, \mathrm{d}x \nonumber\\ [1ex]
& &=\, \frac{\pi}{\cosh \beta - 1}\cdot\left\{ \begin{array}{ll} \big[I_0\big(\sqrt{\beta^2 - 4 \pi^2 v^2}\big) - J_0(2\pi v)\big] &\quad |v| < \frac{\beta}{2\pi}\,,\\ [1ex]
\big[1 - J_0(\beta)\big] &\quad v = \pm \frac{\beta}{2\pi}\,,\\ [1ex]
\big[J_0\big(\sqrt{4 \pi^2 v^2 - \beta^2}\big) - J_0(2\pi v)\big] &\quad |v| > \frac{\beta}{2\pi}\,.
\end{array} \right.
\end{eqnarray}
Note that the Fourier integral in \cite[p.~38, formula 7.57]{Oberh90} reads as follows
\begin{eqnarray*}
\int_0^1  \frac{\cosh\big(\beta \sqrt{1 - x^2}\big)\big)}{\sqrt{1 - x^2}}\,\cos(2\pi v x)\, \mathrm{d}x &=& \int_0^{\pi/2} \cosh(\beta \cos t)\,\cos(2\pi v\,\sin t)\, \mathrm{d}t\\
&=& \frac{\pi}{2}\cdot\left\{ \begin{array}{ll} I_0\big(\sqrt{\beta^2 - 4 \pi^2 v^2}\big) &\quad |v| < \frac{\beta}{2\pi}\,,\\ [1ex]
1 &\quad v = \pm \frac{\beta}{2\pi}\,,\\ [1ex]
J_0\big(\sqrt{4 \pi^2 v^2 - \beta^2}\big)  &\quad |v| > \frac{\beta}{2\pi}\,.
\end{array} \right.
\end{eqnarray*}
Obviously, the unscaled modified $\cosh$-type window function ${\varphi}_{0,\mathrm{mcosh}}:\,\mathbb R \to [0,\,1]$ with the support $[-1,\,1]$ is even and
continuous on $\mathbb R$. Further the restricted function ${\varphi}_{0,\mathrm{mcosh}}|_{[0,\,1]}$ is decreasing. Now we prove that the Fourier transform \eqref{eq:hatomegamcosh(v)} is positive and decreasing in $\big[0,\,\frac{\beta}{2\pi}\big]$. First we remark that by \eqref{eq:hatomegamcosh(v)} it holds
\begin{eqnarray*}
{\hat \varphi}_{0,\mathrm{mcosh}}(0) &=& \frac{\pi}{\cosh \beta - 1}\,I_0(\beta) > 0\,, \\
{\hat \varphi}_{0,\mathrm{mcosh}}\big(\frac{\beta}{2\pi}\big) &=& \frac{2}{\cosh \beta - 1}\,\big[1 - J_0(\beta)\big] > \frac{1}{\cosh \beta - 1} > 0\,,
\end{eqnarray*}
since $J_0(0) = 1$ and $|J_0(\beta)| < \frac{1}{2}$ for $\beta = 2\pi m \big(1 - \frac{1}{2 \sigma}\big) \ge \frac{12\, \pi}{5}$.
Using $I_0'(x) = I_1(x)$ and $J_0'(x) = - J_0(x)$, from \eqref{eq:hatomegamcosh(v)} it follows for all $v \in \big(0,\,\frac{\beta}{2 \pi}\big)$ that
$$
\frac{\mathrm d}{{\mathrm d}v}\,{\hat \varphi}_{0,\mathrm{mcosh}}(v)= \frac{\pi}{\cosh \beta - 1}\Big[- \frac{4 \pi^2 v}{\sqrt{\beta^2 - 4 \pi^2 v^2}} \,I_1\big(\sqrt{\beta^2 - 4 \pi^2 v^2}\big) + 2 \pi\,J_1(2\pi v)\Big]\,.
$$
Since by Lemma \ref{Lemma:inequality1} it holds the inequality
$$
\sqrt{\beta^2 - 4 \pi^2 v^2}\,J_1(2 \pi v) \le 2 \pi v\,I_1\big(\sqrt{\beta^2 - 4 \pi^2 v^2}\big)\,, \quad v  \in \big[0,\,\frac{\beta}{2 \pi}\big]\,,
$$
the Fourier transform \eqref{eq:hatomegamcosh(v)} is decreasing in $\big[0,\,\frac{\beta}{2 \pi}\big]$.

\begin{Lemma}
\label{Lemma:inequality1}
For $x \in [0,\,\beta]$ it holds
\begin{equation}
\label{eq:inequality1}
\sqrt{\beta^2 - x^2}\,J_1(x) \le x\,I_1\big(\sqrt{\beta^2 - x^2}\big)\,.
\end{equation}
\end{Lemma}

\emph{Proof.} Obviously, the inequality \eqref{eq:inequality1} holds for $x = 0$ and $x = \beta$, since $J_1(0) = I_1(0) = 0$. First we prove \eqref{eq:inequality1} for $x \in \big(0,\,\frac{\beta}{\sqrt 2}\big]$.
By the known power expansions
\begin{eqnarray*}
J_1(x) &=& \frac{x}{2}\,\sum_{k=0}^{\infty}\frac{(-1)^k}{4^k\,k!\,(k+1)!}\,x^{2k}\,,\\
I_1(x) &=& \frac{x}{2}\,\sum_{k=0}^{\infty}\frac{1}{4^k\,k!\,(k+1)!}\,x^{2k}\,,
\end{eqnarray*}
we obtain for $x \in \big(0, \,\frac{\beta}{\sqrt 2}\big]$ that
\begin{eqnarray}
\frac{J_1(x)}{x} &=& \frac{1}{2}\,\sum_{k=0}^{\infty}\frac{(-1)^k}{4^k\,k!\,(k+1)!}\,x^{2k}\,, \label{eq:J1(x)/x} \\ [1ex]
\frac{I_1\big(\sqrt{\beta^2 - x^2}\big)}{\sqrt{\beta^2 - x^2}} &=& \frac{1}{2}\,\sum_{k=0}^{\infty}\frac{1}{4^k\,k!\,(k+1)!}\,\big(\beta^2 - x^2\big)^k\,. \label{eq:I1(sqrt)/sqrt}
\end{eqnarray}
Then from $x^2 \le \beta^2 - x^2$ for $x \in \big(0, \,\frac{\beta}{\sqrt 2}\big]$ it follows by \eqref{eq:J1(x)/x} and \eqref{eq:I1(sqrt)/sqrt} that
$$
\frac{J_1(x)}{x} \le \frac{I_1\big(\sqrt{\beta^2 - x^2}\big)}{\sqrt{\beta^2 - x^2}}\,,
$$
i.e., this implies the inequality \eqref{eq:inequality1} for $x \in \big(0, \,\frac{\beta}{\sqrt 2}\big]$.

In the case $x \in \big[\frac{\beta}{\sqrt 2}\,, \,\beta\big)$, we substitute $y:=\sqrt{\beta^2 - x^2} \in \big(0,\,\frac{\beta}{\sqrt 2}\big]$. Then we show
that
$$
\frac{y}{\sqrt{\beta^2 -y^2}}\,  J_1\big(\sqrt{\beta^2 - y^2}\big)  \le I_1(y)\,, \quad y \in \big(0,\,\frac{\beta}{\sqrt 2}\big]\,.
$$
This inequality is fulfilled, since
$$
\big| J_1\big(\sqrt{\beta^2 - y^2}\big) \big| \le \frac{1}{\sqrt 2}
$$
and hence
$$
\frac{y}{\sqrt{\beta^2 -y^2}}\,  J_1\big(\sqrt{\beta^2 - y^2}\big)  \le \frac{y}{\sqrt 2\, \sqrt{\beta^2 - y^2}} \le \frac{y}{\beta} \le I_1(y)\,,
$$
since it holds $\frac{y}{2} \le I_1(y)$ for $y \ge 0$ by $I_1(0) = 0$, $I_1'(0) = \frac{1}{2}$ and $I_1''(y) > 0$.
This completes the proof. \qedsymbol
\medskip

Thus we obtain
$$
{\hat \varphi}_{0,\mathrm{mcosh}}\big(\frac{m}{2 \sigma}\big) = \frac{\pi}{\cosh \beta -1}\,\big[I_0\big(2 \pi m \sqrt{1 - \frac{1}{\sigma}}\big) - J_0\big(\frac{\pi m}{\sigma}\big)\big]\,.
$$
From $m \ge 2$ and $\sigma \in \big[\frac{5}{4},\,2\big]$ it follows that
$$
2 \pi m \sqrt{1 - \frac{1}{\sigma}} \ge 4 \pi \sqrt{1 - \frac{1}{\sigma}} \ge \frac{4 \pi}{\sqrt 5}\,, \quad \frac{\pi m}{\sigma} \ge \pi\,.
$$
Hence it holds
$$
\big|  J_0\big(\frac{\pi m}{\sigma}\big) \big| < \sqrt{\frac{2}{\pi}} \big(\frac{m^2 \pi^2}{\sigma^2} - \frac{1}{4}\big)^{-1/4}\le \sqrt{\frac{2}{\pi}} \big(\pi^2 - \frac{1}{4}\big)^{-1/4} < \frac{1}{2}\,.
$$
This implies that
$$
{\hat \varphi}_{0,\mathrm{mcosh}}\big(\frac{m}{2 \sigma}\big) \ge \frac{\pi}{\cosh \beta -1}\,\big[I_0\big(2 \pi m \sqrt{1 - \frac{1}{\sigma}}\big) - \frac{1}{2}\big]\,.
$$
Further for $v \in \big[m \big(1 - \frac{1}{2 \sigma}\big),\,m \big(1 + \frac{1}{2 \sigma}\big)\big]$ we obtain
$$
|{\hat \varphi}_{0,\mathrm{mcosh}}(v)| = \frac{\pi}{\cosh \beta - 1}\,\big| J_0\big(\sqrt{4\pi^2 v^2 - \beta^2}\big) - J_0(2\pi v)\big| \le c_1 = \frac{3 \pi}{2\,(\cosh \beta -1)}\,,
$$
since it holds
$$
\big| J_0\big(\sqrt{4\pi^2 v^2 - \beta^2}\big)\big| \le 1\,, \quad | J_0(2\pi v)| \le \frac{1}{2}\,.
$$
Thus we can use $c_1 =\frac{3 \pi}{2\,(\cosh \beta - 1)}$ as constant in Theorem \ref{Lemma:esigmaNvarphibounded}.

Finally we determine the decay of \eqref{eq:hatomegamcosh(v)} for $|v| \ge m \big(1 + \frac{1}{2 \sigma}\big)$. We show that
\begin{equation}
\label{eq:est|hatomegamcosh|}
|{\hat \varphi}_{0,\mathrm{mcosh}}(v)| \le \frac{3 \pi m^2 }{2\,(\cosh \beta - 1)}\,\big(1 - \frac{1}{2\sigma}\big)^2\,v^{-2}\,.
\end{equation}
By \eqref{eq:hatomegamcosh(v)} we know that
$$
|{\hat \varphi}_{0,\mathrm{mcosh}}(v)|  = \frac{\pi}{\cosh \beta - 1}\,\big|J_0\big(2\pi v\sqrt{1 - \frac{\beta^2}{4\pi^2v^2}}\big) - J_0(2\pi v)\big|\,.
$$
Using the power expansion of the Bessel function $J_0$, we obtain
$$
J_0\big(2\pi v\sqrt{1 - \frac{\beta^2}{4\pi^2v^2}}\big) - J_0(2\pi v) = \sum_{k=1}^{\infty} \frac{(-1)^k\,(\pi v)^{2k}}{(k!)^2}\,\Big[\big(1 - \frac{\beta^2}{4\pi^2 v^2}\big)^k - 1 \Big]\,.
$$
Since for $k\in \mathbb N$ it holds
$$
0 \le 1 - \big(1 - \frac{\beta^2}{4\pi^2 v^2}\big)^k \le 1 - \big(1 - \frac{\beta^2}{4\pi^2 v^2}\big) = \frac{\beta^2}{4\pi^2 v^2}\,,
$$
we estimate
$$
\big|J_0\big(2\pi v\sqrt{1 - \frac{\beta^2}{4\pi^2v^2}}\big) - J_0(2\pi v)\big| \le |J_0(2 \pi v) - 1|\,\frac{\beta^2}{4\pi^2 v^2}\,.
$$
Since $|J_0(2 \pi v)| \le \frac{1}{2}$ for $v \ge \frac{\beta}{2\pi}$, we receive the inequality \eqref{eq:est|hatomegamcosh|}. Thus
$$
c_2 = \frac{3 \pi m^2 }{2\,(\cosh \beta - 1)}\,\big(1 - \frac{1}{2\sigma}\big)^2\,, \quad \mu = 2
$$
are the corresponding constants of Theorem \ref{Lemma:esigmaNvarphibounded}.

For the modified $\cosh$-type window function \eqref{eq:mcoshwindow}, we estimate the $C(\mathbb T)$-error constant $e_{\sigma}(\varphi_{\mathrm{mcosh}})$.

\begin{Theorem}
\label{Thm:esigma1(varphimcosh)}
Let $\sigma \in \big[\frac{5}{4},\,2\big]$, $m \in {\mathbb N}\setminus \{1\}$, and $N_1 = \sigma_1 N\in 2 \mathbb N$ with $2 m \ll N_1$.

Then the $C(\mathbb T)$-error constant $e_{\sigma}(\varphi_{\mathrm{mcosh}})$ has the upper bound
\begin{equation}
\label{eq:esigma1(varphimcosh)}
e_{\sigma}(\varphi_{\mathrm{mcosh}}) \le \frac{21}{4}\,\big[ I_0\big(2 \pi m\,\sqrt{1 - \frac{1}{\sigma}}\big) + \frac{1}{2}\big]^{-1}\,.
\end{equation}
\end{Theorem}

\emph{Proof.} From Theorem \ref{Lemma:esigmaNvarphibounded} with $\mu = 2$ it follows that
$$
e_{\sigma}(\varphi_{\mathrm{mcosh}}) \le \frac{1}{{\hat \varphi}_{0,\mathrm{mcosh}}\big(\frac{m}{2 \sigma}\big)}\,\big[2 c_1 + \frac{2c_2}{m^2}\,\big(1 - \frac{1}{2 \sigma}\big)^{-1}\big]\,.
$$
Since
$$
{\hat \varphi}_{0,\mathrm{mcosh}}\big(\frac{m}{2 \sigma}\big) \ge \frac{\pi}{\cosh \beta -1}\,\big[I_0\big(2 \pi m \sqrt{1 - \frac{1}{\sigma}}\big) - \frac{1}{2}\big]
$$
and
$$
c_1 = \frac{3 \pi}{2\,(\cosh \beta -1)}\,, \quad c_2 = \frac{3 \pi m^2 }{2\,(\cosh \beta - 1)}\,\big(1 - \frac{1}{2\sigma}\big)^2\,,
$$
we obtain the estimate \eqref{eq:esigma1(varphimcosh)}. Note that by $\sigma \in \big[\frac{5}{4},\,2\big]$ it holds
$$
3 + 3\,\big(1 - \frac{1}{2 \sigma}\big) \le \frac{21}{4}\,.
$$
By \cite{YaChu16} the modified Bessel function $I_0(x)$ fulfills the inequalities
$$
\frac{{\mathrm e}^x}{1 + 2 x} < I_0(x) < \frac{{\mathrm e}^x}{\sqrt{1 + 2 x}}\,, \quad x > 0\,.
$$
If we apply this result, then we obtain
$$
e_{\sigma}(\varphi_{\mathrm{mcosh}}) \le \frac{21}{4}\,\big(1 + 4\pi m \sqrt{1- \frac{1}{\sigma}}\big)\Big[{\mathrm e}^{2 \pi m\,\sqrt{1 - 1/\sigma}} + \frac{1}{2} + 2 \pi m \sqrt{1 - \frac{1}{\sigma}}\Big]^{-1}\,.
$$
This completes the proof. \qedsymbol

\medskip

Finally we consider in our numerical examples also  the \emph{modified} $\mathrm{exp}$-\emph{type window function}
\begin{equation*}
\varphi_{\mathrm{mexp}}(x) := \varphi_{0,\mathrm{mexp}}\big(\frac{N_1 x}{m}\big)\,, \quad x \in \mathbb R\,,
\end{equation*}
where it holds
\begin{equation*} 
\varphi_{0,\mathrm{mexp}} (x) := \left\{ \begin{array}{ll}
\frac{1}{\exp \beta - 1}\,\frac{\exp \big(\beta \sqrt{1 - x^2}\big)- 1}{\sqrt{1 - x^2}}    &\quad  x \in ( - 1,\, 1)\,,\\
0                   & \quad x \in \mathbb R \setminus ( - 1,\, 1)
\end{array} \right.
\end{equation*}
and the \emph{modified} $\sinh$-\emph{type window function}
\begin{equation*}
\varphi_{\mathrm{msinh}}(x) := \varphi_{0,\mathrm{msinh}}\big(\frac{N_1 x}{m}\big)\,, \quad x \in \mathbb R\,,
\end{equation*}
where it holds
\begin{equation*}
\varphi_{0,\mathrm{msinh}} (x) := \left\{ \begin{array}{ll}
\frac{1}{\sinh \beta}\,\frac{\sinh \big(\beta \sqrt{1 - x^2}\big)}{\sqrt{1 - x^2}}    &\quad  x \in ( - 1,\, 1)\,,\\
0                   & \quad x \in \mathbb R \setminus ( - 1,\, 1).
\end{array} \right.
\end{equation*}

Note that the main idea to consider these modified window functions, comes from the Fourier transform of the Kaiser--Besser window, see \cite[Remark 1]{post02}. The modified $\sinh$-type window function $\varphi_{\mathrm{msinh}}$ is used in the NFFT package \cite{nfft3}, and gives very good error bounds.
}

\section*{Conclusion}

In this paper, we present explicit error estimates for the  NFFT  with continuous, compactly supported window function $\varphi \in \Phi_{m,N_1}$. Such window functions simplify the algorithm for NFFT,
since the truncation error of the NFFT vanishes. Using the $C(\mathbb T)$-error constant of $\varphi$, we propose a unified approach to error estimates of the NFFT with nonequispaced spatial data and equispaced frequencies
as well as of the NFFT with nonequispaced frequencies and equispaced spatial data. Further we discuss the connection with a modified Paley--Wiener theorem.

We present two techniques to find upper bounds of the $C(\mathbb T)$-error constant. The second method which uses the scaling structure
of the window function $\varphi$, shows that the constant $e_{\sigma,N}(\varphi)$ is bounded for all $N \in 2 \mathbb N$. We see that $e_{\sigma}(\varphi)$ depends essentially {  on} the decay of the Fourier transform
${\hat \varphi}(v)$ for $|v| \to \infty$ and the positive size of ${\hat \varphi}(v)$ for small frequencies. For the (modified) B-spline, algebraic, Bessel, $\sinh$-type, and modified $\cosh$-type window functions, we present
new explicit upper bounds of the corresponding $C(\mathbb T)$-error constants. Here we use the fact that the Fourier transforms of these window functions are explicitly known. It is remarkable that the
$C(\mathbb T)$-error constants of Bessel, $\sinh$-type, and modified $\cosh$-type window function decay exponentially with respect to $m$. Numerical experiments verify the different behavior of the $C(\mathbb T)$-error constants for these window functions, see Figure \ref{Fig.Vergl}. We point out that the modified $\cosh$-type, $\exp$-type, and $\sinh$-type window functions give the best error constants, see Figure \ref{Fig.Vergl2}. The modified $\sinh$-type window function, see Section \ref{Sec:modwindow}, is used in the NFFT package \cite{nfft3}.


\begin{figure}[ht]
	\begin{center}
		\begin{tikzpicture}[scale=0.5]
		\begin{axis}[
		ymode=log,
		title={Error for $\sigma=1.25$},
		xlabel={$m$},
		xmin=1.8, xmax=6.2,
		ymin=0.8*10^(-6), ymax=5*10^(-1),
		xtick={2,3,4,5,6},
		ytick={10^(-1),10^(-2),10^(-3),10^(-4), 10^(-5),10^(-6)},
		yticklabels={$10^{-1}$, $10^{-2}$, $10^{-3}$, $10^{-4}$, $10^{-5}$, $ 10^{-6}$},
		legend pos=south west,
		ymajorgrids=true,
		grid style=dashed,
		]
		\addplot[color=blue, mark=square]
		coordinates {(2, 2.1004e-01) (3,8.8620e-02) (4,3.9079e-02)  (5, 1.7346e-02 )  (6,7.7077e-03)};
		\addplot[color=red, mark=o]
		coordinates {(2,2.3270e-01)  (3,8.3879e-02 )   (4,2.6817e-02) (5,7.3679e-03 ) (6,1.5376e-03)};
		\addplot[color=green, mark=*]
		coordinates {(2,1.5318e-01)  (3,2.5431e-02  )   (4, 3.2884e-03 ) (5,3.6505e-04  ) (6,3.6565e-05)};
		\addplot[color=black, mark=triangle]
		coordinates {(2, 6.6976e-02)  (3,7.1278e-03)   (4,6.5055e-04)  (5,5.4320e-05) (6,4.2561e-06)};
		\addplot[color=cyan, mark=otimes]
		coordinates {(2, 1.687e-02)  (3,2.1601e-03)   (4,2.9749e-04)  (5,1.3047e-05) (6,1.2548e-06)};
		\legend{$\varphi_{\mathrm{B}}$, $\varphi_{\mathrm{alg}}$,$\varphi_{\mathrm{Bessel}}$, $\varphi_{\sinh}$, $\varphi_{\mathrm{mcosh}}$}
		\end{axis}
		\end{tikzpicture}
		\begin{tikzpicture}[scale=0.5]
		\begin{axis}[
		ymode=log,
		title={Error for $\sigma=1.5$},
		xlabel={$m$},
		xmin=1.8, xmax=6.2,
		ymin=0.8*10^(-8), ymax=5*10^(-2),
		xtick={2,3,4,5,6},
		ytick={10^(-2),10^(-3),10^(-4),10^(-5), 10^(-6), 10^(-7), 10^(-8)},
		yticklabels={$10^{-2}$, $10^{-3}$, $10^{-4}$, $10^{-5}$, $10^{-6}$, $10^{-7}$, $10^{-8}$},
		legend pos=south west,
		ymajorgrids=true,
		grid style=dashed,
		]
		\addplot[color=blue, mark=square]
		coordinates {(2, 6.8961e-02) (3,1.5948e-02) (4,3.9243e-03 )  (5, 9.7762e-04 )  (6,2.4420e-04)};
		\addplot[color=red, mark=o]
		coordinates {(2,7.7777e-02)  (3,9.7487e-03 )   (4,9.3175e-04) (5,5.0595e-04 ) (6,1.7182e-04)};
		\addplot[color=green, mark=*]
		coordinates {(2, 5.6157e-02 )  (3,4.3081e-03  )   (4, 2.5194e-04 ) (5,1.2515e-05   ) (6,5.5830e-07)};
		\addplot[color=black, mark=triangle]
		coordinates {(2, 1.9459e-02)  (3,8.9016e-04)   (4,3.5864e-05)  (5,1.3237e-06) (6,4.5903e-08)};
		\addplot[color=cyan, mark=otimes]
		coordinates {(2, 9.3633e-03)  (3,3.0073e-04)   (4,1.0503e-05)  (5,2.5092e-07) (6,9.3504e-09)};
		\legend{$\varphi_{\mathrm{B}}$, $\varphi_{\mathrm{alg}}$,$\varphi_{\mathrm{Bessel}}$, $\varphi_{\sinh}$,
			$\varphi_{\mathrm{mcosh}}$}
		\end{axis}
		\end{tikzpicture}
		\begin{tikzpicture}[scale=0.5]
		\begin{axis}[
		ymode=log,
		title={Error for $\sigma=2$},
		xlabel={$m$},
		xmin=1.8, xmax=6.2,
		ymin=0.5*10^(-10), ymax=2*10^(-2),
		xtick={2,3,4,5,6},
		ytick={10^(-2), 10^(-3), 10^(-4), 10^(-5),10^(-6),10^(-7), 10^(-8), 10^(-9),10^(-10)},
		yticklabels={$10^{-2}$, $10^{-3}$, $10^{-4}$,$10^{-5}$,$10^{-6}$,$10^{-7}$, $10^{-8}$,$10^{-9}$,$10^{-10}$},
		legend pos=south west,
		ymajorgrids=true,
		grid style=dashed,
		]
		\addplot[color=blue, mark=square]
		coordinates {(2, 1.4678e-02 ) (3,1.4471e-03 ) (4,1.5518e-04)  (5, 1.7041e-05  )  (6,1.8859e-06)};
		\addplot[color=red, mark=o]
		coordinates {(2,1.1934e-02)  (3,2.0136e-03 )   (4,4.4476e-04 ) (5,2.9620e-05 ) (6,3.4013e-06)};
		\addplot[color=green, mark=*]
		coordinates {(2, 1.8493e-02  )  (3,6.4916e-04  )   (4, 1.7039e-05  ) (5,3.7777e-07  ) (6,7.5002e-09)};
		\addplot[color=black, mark=triangle]
		coordinates {(2,5.1243e-03) (3,1.0287e-04 )  (4,1.8467e-06) (5,3.0197e-08) (6,4.6553e-10)};
		\addplot[color=cyan, mark=otimes]
		coordinates {(2,3.5681e-03) (3,2.6533e-05 )  (4,5.1780e-07) (5,4.0321e-09) (6,1.1991e-10)};
		\legend{$\varphi_{\mathrm{B}}$, $\varphi_{\mathrm{alg}}$,$\varphi_{\mathrm{Bessel}}$, $\varphi_{\sinh}$,
			$\varphi_{\mathrm{mcosh}}$}
		\end{axis}
		\end{tikzpicture}
		
		\caption{The constants $e_{\sigma,1024}({\varphi})$ of the different window functions  with shape pa\-rameter $\beta=\pi m(2-1/\sigma)$
			for  $\sigma\in \{1.25, 1.5,2 \}$ and  $m\in \{2,3,4,5,6\}$.}
		\label{Fig.Vergl}
	\end{center}	
\end{figure}
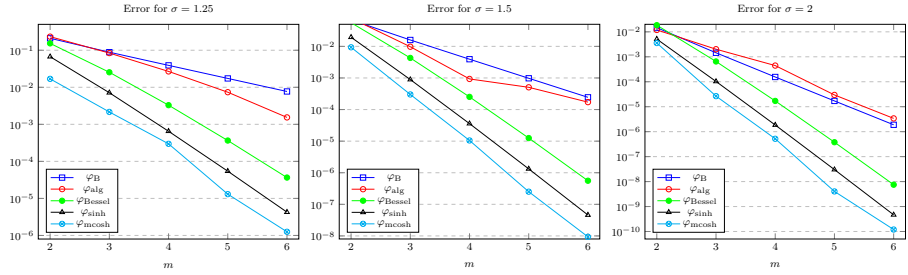

\begin{figure}[ht]
	\begin{center}
		\begin{tikzpicture}[scale=0.5]
		\begin{axis}[
		domain=2:6,
		ymode=log,
		title={Error for $\sigma=1.25$},
		xlabel={$m$},
		xmin=1.8, xmax=6.2,
		ymin=0.8*10^(-6), ymax=5*10^(-1),
		xtick={2,3,4,5,6},
		ytick={10^(-1),10^(-2),10^(-3),10^(-4), 10^(-5),10^(-6)},
		yticklabels={$10^{-1}$, $10^{-2}$, $10^{-3}$, $10^{-4}$, $10^{-5}$, $ 10^{-6}$},
		legend pos=south west,
		ymajorgrids=true,
		grid style=dashed,
		]
		\addplot[color=black, mark=*]
		coordinates {(2, 6.6976e-02)  (3,7.1278e-03) (4,6.5055e-04) (5,5.4320e-05) (6,4.2561e-06)};
		\addplot[color=cyan, mark=*]
		coordinates {(2, 2.76e-02) (3,1.9713e-03) (4,2.35e-04)  (5, 9.0427e-06)  (6,9.86e-07)};
		 \addplot[color=black,mark=triangle] {(24*sqrt(x*x)+3)*exp(-2*pi*x*sqrt(1-1/1.25))};	
		\addplot[color=cyan, mark=otimes]
		{21/4/(  (exp(2*pi*x*sqrt(1-1/1.25))/sqrt(1+2*(2*pi*x*sqrt(1-1/1.25))))+1/2)};	
		\legend{
			\Large{$\varphi_{\sinh}$=$\varphi_{\cosh}$=$\varphi_{\mathrm{exp}}$},
			\Large {$\varphi_{\mathrm{msinh}}$=$\varphi_{\mathrm{mcosh}}$=$\varphi_{\mathrm{mexp}}$},
			estimate Theorem \ref{Thm:esigmaNsinh},
			estimate Theorem \ref{Thm:esigma1(varphimcosh)}
		}
		\end{axis}
		\end{tikzpicture}
		\begin{tikzpicture}[scale=0.5]
		\begin{axis}[
		domain=2:6,
		ymode=log,
		title={Error for $\sigma=1.5$},
		xlabel={$m$},
		xmin=1.8, xmax=6.2,
		ymin=0.8*10^(-8), ymax=5*10^(-2),
		xtick={2,3,4,5,6},
		ytick={10^(-2),10^(-3),10^(-4),10^(-5), 10^(-6), 10^(-7), 10^(-8)},
		yticklabels={$10^{-2}$, $10^{-3}$, $10^{-4}$, $10^{-5}$, $10^{-6}$, $10^{-7}$, $10^{-8}$},
		legend pos=south west,
		ymajorgrids=true,
		grid style=dashed,
		]
		\addplot[color=black, mark=*]
		coordinates {(2,1.9433e-02 ) (3,8.9016e-04) (4,3.5865e-05) (5,1.3237e-06) (6,4.5903e-08)};
		\addplot[color=cyan, mark=*]
		coordinates {(2,9.49e-03) (3,2.007e-04) (4,9.35e-06 ) (5,2.66e-07) (6,9.3504e-09)};
		 \addplot[color=black,mark=triangle] {(24*sqrt(x*x)+3)*exp(-2*pi*x*sqrt(1-1/1.5))};	
		\addplot[color=cyan, mark=otimes]
		{21/4/(  (exp(2*pi*x*sqrt(1-1/1.5))/sqrt(1+2*(2*pi*x*sqrt(1-1/1.5))))+1/2)};
		\legend{
			\Large{$\varphi_{\sinh}$=$\varphi_{\cosh}$=$\varphi_{\mathrm{exp}}$},
			\Large {$\varphi_{\mathrm{msinh}}$=$\varphi_{\mathrm{mcosh}}$=$\varphi_{\mathrm{mexp}}$},
			estimate Theorem \ref{Thm:esigmaNsinh},
			estimate Theorem \ref{Thm:esigma1(varphimcosh)}
		}
		\end{axis}
		\end{tikzpicture}
		\begin{tikzpicture}[scale=0.5]
		\begin{axis}[
		domain=2:6,
		ymode=log,
		title={Error for $\sigma=2$},
		xlabel={$m$},
		xmin=1.8, xmax=6.2,
		ymin=0.5*10^(-10), ymax=2*10^(-2),
		xtick={2,3,4,5,6},
		ytick={10^(-2), 10^(-3), 10^(-4), 10^(-5),10^(-6),10^(-7), 10^(-8), 10^(-9),10^(-10)},
		yticklabels={$10^{-2}$, $10^{-3}$, $10^{-4}$,$10^{-5}$,$10^{-6}$,$10^{-7}$, $10^{-8}$,$10^{-9}$,$10^{-10}$},
		legend pos=south west,
		ymajorgrids=true,
		grid style=dashed,
		]
		\addplot[color=black, mark=*]
		coordinates {(2,5.1e-03) (3,1.0287e-04 )  (4,1.8467e-06) (5,3.0197e-08) (6,4.6553e-10)};
		\addplot[color=cyan, mark=*]
		coordinates {(2, 3.09e-03)  (3,2.6703e-05)  (4,3.282e-07) (5,4.867e-09) (6, 8.482e-11)};
	    \addplot[color=black,mark=triangle] {(24*sqrt(x*x)+3)*exp(-2*pi*x*sqrt(1-1/2))};	
	 \addplot[color=cyan, mark=otimes] 
	{21/4/(  (exp(2*pi*x*sqrt(1-1/2))/sqrt(1+2*(2*pi*x*sqrt(1-1/2))))+1/2)};		
		\legend{
			\Large{$\varphi_{\sinh}$=$\varphi_{\cosh}$=$\varphi_{\mathrm{exp}}$},
			\Large {$\varphi_{\mathrm{msinh}}$=$\varphi_{\mathrm{mcosh}}$=$\varphi_{\mathrm{mexp}}$},
			estimate Theorem \ref{Thm:esigmaNsinh},
			estimate Theorem \ref{Thm:esigma1(varphimcosh)}
		}
		\end{axis}
		\end{tikzpicture}
		
		\caption{The constants $e_{\sigma,1024}({\varphi})$ of the different window functions  with shape pa\-rameter $\beta=\pi m(2-1/\sigma)$
			for  $\sigma\in \{1.25, 1.5,2 \}$ and  $m\in \{2,\,3,\,4,\,5,\,6\}$ as well as the upper bounds of the $C(\mathbb T)$-error constants $e_{\sigma}(\varphi)$.}
		\label{Fig.Vergl2}
	\end{center}
\end{figure}
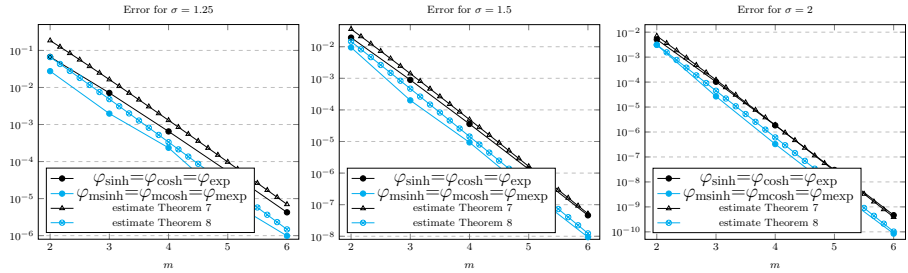

\section*{Acknowledgments}

The authors would like to thank  J\"{u}rgen Rossmann and Winfried Sickel  for several fruitful discussions on the topic. {  Further, the authors would like to thank the reviewers for constructive advices to improve the representation of the paper.}
The first author acknowledges funding  by Deutsche Forschungsgemeinschaft (German Research Foundation) -- Project--ID 416228727 -- SFB 1410.

\end{document}